\documentclass[10pt, a4paper, reqno]{amsart}
\usepackage{style}

\begin{document}

\begin{abstract}
    We identify a class of singular algebraic foliations whose leaves through singular points retain regularity.
    The proof consists in showing existence of residual gerbes for certain formal stacks, which do not enjoy smooth presentations.
    As applications, we extend a theorem of Cerveau to the case where the ambient scheme is not smooth and we give a proof of the Zariski--Lipman conjecture for varieties with terminal singularities, which does not rely on existence of resolution of singularities.
\end{abstract}

\maketitle

\tableofcontents
\section*{Introduction}
\label{sec:introduction}

If a modern algebraic geometer were forced to choose between working with smoothness hypotheses or compactness hypotheses, they would probably choose the latter.
Indeed, it has long been realised that allowing singularities to come into the picture is manageable, and sometimes even desirable, so long as they are controlled.
Foliations in Algebraic Geometry are no exception: it is extremely rare to find a smooth algebraic foliation on a projective variety.
In the literature, the approaches to the study and control of foliation singularities fall into two categories:
\begin{enumerate}
    \item via singularities of the Minimal Model Program, as pioneered by McQuillan in \cite{MR2435846}, and
    \item via formal groupoids, as more recently introduced by McQuillan in \cite{MR4507257}.
\end{enumerate}
The first approach is the direct analogue of the largely successful birational geometric techniques of the Minimal Model Program.
The second approach should perhaps be thought as the analogue of studying varieties with quotient singularities.
We do not try to conceal the fact that we prefer the second approach.
In order for this to work, we limit ourselves to considering foliations whose tangent sheaf $\mathscr{F} \subseteq \mathscr{T}_X$ is closed under the Lie bracket and \emph{locally free}, although perhaps unsaturated, on an ambient scheme $X$ with arbitrary singularities or non-reducedness.
As an example, we note that any vector field, whether singular or non-singular, belongs to this category.
We refer the reader to Example \ref{ex:foliations} for further instances of locally free foliations and to Remark \ref{rem:foliation_resolution} for a plausible potential treatment of the general non-locally free case.

\subsection*{Main Results}
When $X$ is a scheme over an algebraically closed field of characteristic zero and $\mathscr{F}$ is a locally free foliation, we show the following results.

\begin{corollary*}[\ref{cor:smoothness_leaf}]
    The formal leaf $L$ through any closed point $x$, defined as the minimal $\mathscr{F}$-invariant subscheme through $x$, exists and is smooth.
\end{corollary*}

The following corollary is an extension of a theorem of Cerveau (\cite[Theorem 1.1]{MR552968}) to the case where $X$ is a singular scheme.
We thank Jorge Vit\'{o}rio Pereira for introducing us to this paper.

\begin{corollary*}[\ref{cor:cerveau_foliation}]
    There exists a transversal $V$ to the formal leaf $L$ through $x$, such that $X$ locally decomposes as a product $L \times V$, which is compatible with the foliated structure.
\end{corollary*}

\begin{lemma*}[\ref{lem:blow_up_foliation}]
    Blowing up with respect to $\mathscr{F}$-invariant subschemes preserves the foliated structure.
\end{lemma*}

As a consequence, we show some applications. \par

Firstly, we generalise a theorem of Hart (\cite[Corollary 2]{MR349654}) to the relative setting.

\begin{lemma*}[\ref{lem:fitting_ideals_connection}]
    Given a morphism of schemes $f : X \rightarrow Y$, any derivation $\partial$ of $X$ over $Y$ preserves the Fitting ideals of $\Omega_{X/Y}^1$.
\end{lemma*}

Secondly, we generalise a well-known characterisation of terminal foliations by curves to higher rank.

\begin{proposition*}[\ref{prop:terminal_singularities}]
    If a locally free saturated foliation $\mathscr{F}$ on a normal variety $X$ has at worst terminal singularities, then both $X$ and $\mathscr{F}$ are smooth.
\end{proposition*}

Finally we give a proof of the Zariski--Lipman conjecture for terminal varieties.
This does not generalise anything.
In fact the statement is known for log-canonical singularities (\cite[Theorem 5.2]{MR3239620} or \cite[Corollary 1.3]{MR3247804}).
It is nonetheless here included as the proof does not depend on the existence of resolution of singularities.

\begin{proposition*}[\ref{prop:zariski_lipman}]
    If the tangent sheaf of a normal terminal variety $X$ is locally free, then $X$ is smooth.
\end{proposition*}

In order to avoid misunderstandings, we urge the reader to consult the attached precise statements, which include finiteness and Noetherian hypotheses. \par

\subsection*{Infinitesimal Stacks}

Most of the above results are applications of Theorem \ref{thm:minimal_presentation}, a result concerning existence of residual gerbes and minimal presentations of infinitesimal stacks defined over a general base scheme.
We then translate our results in characteristic zero, where there is an equivalence of categories between infinitesimal stacks and foliations.
In fact, we take the opportunity to display a text converter between the two categories:

\begin{center}
    \renewcommand{\arraystretch}{1.5}
    \begin{tabular}{
        |>{\centering\arraybackslash}m{0.32\textwidth}|
         >{\centering\arraybackslash}m{0.26\textwidth}|
         >{\centering\arraybackslash}m{0.32\textwidth}|}
    \hline
    \textbf{Infinitesimal Stacks} & \textbf{$\longleftrightarrow$} & \textbf{Foliations} \\
    \hline\hline
    Formally smooth infinitesimal groupoids $\mathcal{X} = [R \rightrightarrows X]$ & Proposition \ref{prop:foliation_groupoid} & Locally free foliations $\mathscr{F}$ on $X$ \\
    \hline
    Formal $\hat{\mathbb{G}}_m$-actions on $X$ & Example \ref{ex:foliations} & Vector fields $\partial$ on $X$ \\
    \hline
    Stabiliser group $G \rightarrow X$ & Lemma \ref{lem:stabiliser_surjection} & Cokernel of $\Omega_X^1 \rightarrow \mathscr{F}^{\vee}$ \\
    \hline
    Residual gerbe of a point $x \in \mathcal{X}$ & Remark \ref{rem:residual_gerbes} & Leaf $L$ through $x \in  X$ \\
    \hline
    Minimal presentations of $\mathcal{X}$ at $x$ & Theorem \ref{thm:minimal_presentation} & Transversals $V$ of $\mathscr{F}$ through $x$ \\
    \hline
    Locally closed substacks & Lemma \ref{lem:invariance_equivalence} & $\mathscr{F}$-invariant subschemes \\
    \hline
    Representable morphisms of stacks $f : \mathcal{W} \rightarrow \mathcal{X}$ & Lemma \ref{lem:functorial_foliations} & Morphisms of foliated spaces $g : (W, \mathscr{G}) \rightarrow (X, \mathscr{F})$ such that $\mathscr{G} = f^*\mathscr{F}$ \\
    \hline
    Coherent sheaves $\mathscr{E}$ on $\mathcal{X}$ & Proposition \ref{prop:equivariance_equivalence} & Coherent sheaves $\mathscr{E}$ on $X$ with a flat partial connection in $\mathscr{F}$ \\
    \hline
    Cotangent complex of $\mathcal{X}$ & Example \ref{ex:bott_connection} & Bott connection on the normal bundle $\mathscr{N}_{\mathscr{F}}$ \\
    \hline
    \end{tabular}
\end{center}

\subsection*{Main Idea}

We give a crude idea of the proof.
We first point out that the proof of existence of residual gerbes for algebraic stacks cannot be adapted to our infinitesimal setting.
Indeed, in the case of algebraic stacks, a crucial step is to use generic flatness of the morphism between $x$ and the stack-theoretic image in $\mathcal{X}$.
However, formally smooth presentations of infinitesimal stacks are almost never of finite presentation, and generic flatness may not hold. \par

Instead, our proof uses flattening stratifications in non-finitely presented contexts.
We apply the following principle: given an $R$-equivariant morphism $g : W \rightarrow X$, the non-flat locus of $g$ in $X$ should be $R$-invariant.
When $X$ is the leaf of a foliation, by definition, there cannot be any $R$-invariant subschemes, thus $g$ is flat.
Our proof also works in the case of algebraic stacks. \par

In fact, most of the arguments in this paper involve considering appropriate functorial constructions, which yield $R$-invariant subschemes.
We then try to exploit these special subschemes.

\subsection*{Outline}

In \S \ref{sec:geometry_formal_groupoids}, we discuss some preliminary results concerning coherent sheaves on formal groupoids, invariance of subschemes and equivariance of morphisms.
In \S \ref{sec:functorial_constructions}, we carry out the required functorial constructions: reduction, normalisation, Fitting ideals, flattening stratifications, blowing up and functorial resolutions.
In \S \ref{sec:formal_leaves}, we define formal leaves and apply the previous constructions to deduce regularity and flatness properties.
In \S \ref{sec:minimal_presentations}, we define transversals and we prove the main theorem.
In \S \ref{sec:foliations}, we apply our methods to deduce the aforementioned results about foliations over a field of characteristic zero.

\subsection*{Assumptions and Notation}

Throughout the article, we work in the category of formal schemes as defined in \cite{bongiorno1}.
In particular, a formal scheme is a locally topologically ringed space which is locally isomorphic to the formal spectrum of an adic ring with a finitely generated ideal of definition. \par

We let $S$ denote a base formal scheme which, for simplicity, may be thought as the spectrum of a field, and we let $R \rightrightarrows X$ denote a \emph{formal groupoid} over $S$, i.e. a groupoid object in the category of formal schemes over $S$.
The groupoid structure entails
\begin{enumerate}
    \item[($s$)] a source morphism $s : R \rightarrow X$,
    \item[($t$)] a target morphism $t : R \rightarrow X$,
    \item[($e$)] a unit morphism $e : X \rightarrow R$,
    \item[($i$)] an inverse morphism $i : R \rightarrow R$, and
    \item[($c$)] a composition morphism $c : R \times_{(s,t)} R \rightarrow R$
\end{enumerate}
satisfying various compatibility conditions.
We also define $j := t \times s : R \rightarrow X \times_S X$ to be the diagonal morphism.
An \emph{infinitesimal groupoid} is a formal groupoid whose unit morphism is a thickening of formal schemes, or equivalently, whose source morphism induces a homeomorphism of topological spaces.
We refer the reader to \cite[\S 1]{bongiorno2} for further details. \par

We will be mostly concerned with locally Noetherian groupoids.
\begin{situation-g}
\phantomsection\label{sit:global}
    $X$ is a formal scheme locally of formal finite presentation over a Noetherian formal scheme $S$ and $R$ is an infinitesimal groupoid on $X$ over $S$ whose source morphism is formally smooth and locally of formal finite presentation.
\end{situation-g}

We refer the reader to \cite[\S 4]{bongiorno1} for the definition and properties of morphisms locally of formal finite presentation. \par

We will often pick a closed point $x \in X$ and study the formal neighbourhood around it.
This is permitted since, in Situation \hyperref[sit:global]{Global}, it follows that both $X$ and $R$ are locally Noetherian formal schemes (\cite[Lemma \ref*{inf-lem:noetherian_ascent}]{bongiorno1}) and both closed immersions $x \in X$ and $e(x) \in R$ are locally of formal finite presentation, hence we may consider their respective infinitesimal neighbourhoods $\hat{X}$ and $\hat{R}$.
In fact, $\hat{R}$ is naturally isomorphic to the restriction of the groupoid $R$ via the closed immersion $x \in X$ (\cite[Lemma \ref*{frob-lem:formal_restriction}]{bongiorno2}), in particular $\hat{R}$ is an infinitesimal groupoid.
Furthermore, the induced source morphism $\hat{s} : \hat{R} \rightarrow \hat{X}$ is formally smooth and locally of formal finite presentation (see \cite[Construction \ref*{frob-cons:infinitesimal_from_formal}]{bongiorno2}), thus we reduce to

\begin{situation-l}
\phantomsection\label{sit:local}
    $X$ is a formal scheme over a formal scheme $S$ and the topological space of $X$ consists of a single closed point $x$.
    $R$ is an infinitesimal groupoid on $X$ over $S$ whose source morphism is formally smooth and locally of formal finite presentation.
\end{situation-l}

In this case, it follows that both $X = \spf A$ and $R = \spf \Gamma$ are affine formal schemes (\cite[Lemma \ref*{inf-lem:image_thickening}]{bongiorno1}), and the closed immersion $x \in X$ corresponds to the maximal adic ideal of definition $\mathfrak{m}$ of $A$.

We will be using, without further mention, the following facts:
\begin{enumerate}
    \item if $A$ is a Noetherian local ring, its adic completion $A \rightarrow \hat{A}$ is faithfully flat,
    \item if $\varphi : A \rightarrow B$ is flat local morphism of local rings, it is faithfully flat (\cite[\href{https://stacks.math.columbia.edu/tag/00HR}{Lemma 00HR}]{stacks-project}),
    \item if $\varphi : A \rightarrow B$ is local morphism of adic Noetherian local rings, then $\varphi$ is formally smooth if and only if $\varphi$ is faithfully flat and its central fibre is formally smooth (\cite[\href{https://stacks.math.columbia.edu/tag/07NQ}{Proposition 07NQ}]{stacks-project}),
    \item if $f : X \rightarrow Y$ is a formally smooth morphism locally of formal finite presentation of locally Noetherian formal schemes, then $f$ is flat and the relative sheaf of differentials $\Omega_{X/Y}^1$ is a locally free $\mathscr{O}_X$-module of finite rank (\cite[Proposition 4.8]{MR2313672}),
    \item if $f : X \rightarrow Y$ is a thickening of formal schemes and $X$ is affine, so is $Y$ (\cite[Lemma \ref*{inf-lem:image_thickening}]{bongiorno1}), and
    \item if $f : X \rightarrow Y$ is a morphism locally of formal finite presentation of formal schemes and $Y$ is locally Noetherian, so is $X$ (\cite[Lemma \ref*{inf-lem:noetherian_ascent}]{bongiorno1}).
\end{enumerate}

\subsection*{Acknowledgements}

I would like to thank my doctoral advisor Paolo Cascini for the inspiring course given at Tsinghua Univeristy.
I would also like to thank Jarod Alper, Marc Besson, Caucher Birkar, Alessio Bottini, Riccardo Carini, Priyankur Chaudhuri, St\'{e}phane Druel, Gabriel Fazoli, Benjamin Hennion, Crislaine Kuster, Li Shang, Theodoros Papazachariou, Jorge Vit\'{o}rio Pereira, Quentin Posva, Sheng Mao and Stefania Vassiliadis for helpful advice and discussions.
I thank my home institution of Tsinghua University where most of the research work was carried out, and Instituto de Matem\'{a}tica Pura e Aplicada for their summer postdoctoral program.
I acknowledge funding from both institutions.

\section{Geometry on Formal Groupoids}
\label{sec:geometry_formal_groupoids}

In this section, we collect some facts about equivariant sheaves on formal schemes with respect to an infinitesimal groupoid.
We first define coherent sheaves on formal schemes in \S \ref{subsec:coherent_sheaves} and then proceed with defining equivariant structures in \S \ref{subsec:equivariant_structures}.
We provide several examples of such structures.

\subsection{Coherent Sheaves}
\label{subsec:coherent_sheaves}

We follow \cite[\href{https://stacks.math.columbia.edu/tag/0EHN}{Section 0EHN}]{stacks-project} for basic notions of coherent sheaves on affine Noetherian formal schemes.

\begin{definition}
\label{def:adic_modules}
    Let $A$ be an adic Noetherian ring and let $I$ be an ideal of definition of $A$.
    An \emph{adic $A$-module} $M$ is a finitely presented topological $A$-module such that
     \begin{enumerate}
        \item $\{I^{n+1} M \}_{n \in \mathbb{N}}$ is a fundamental system of submodules for the topology of $M$, and
        \item $M$ is a complete and separated topological space.
    \end{enumerate}
    A \emph{morphism of adic modules} is a continuous morphism of topological modules.
\end{definition}

As usual, being adic does not depend on the choice of ideal of definition.
We may define a tensor product operation $\_ \hat{\otimes}_A \_$ of two adic $A$-modules by first performing usual tensor product over $A$ and then completing with respect to the induced topology.
The resulting module is an adic $A$-module.

\begin{definition}
\label{def:coherent_sheaves}
    Let $X$ be a locally Noetherian formal scheme.
    A \emph{coherent sheaf} $\mathscr{F}$ on $X$ is a sheaf of topological modules over $\mathscr{O}_X$ such that for all affine open sets $U \subseteq X$, $\mathscr{F}(U)$ is an adic $\mathscr{O}_X(U)$-module.
    A \emph{morphism of coherent sheaves} is a continuous morphism of sheaves of topological modules.
\end{definition}

It is enough to check coherence of a module on an affine open cover.
We may define tensor product of coherent sheaves by first defining it locally as above and then glueing.

\begin{remark}
\label{rem:coherent_sheaves}
    Let $X$ be a Noetherian formal scheme with global ideal sheaf of definition $\mathscr{I} \subseteq \mathscr{O}_X$.
    Then the datum of a coherent sheaf $\mathscr{F}$ on $X$ is equivalent to the datum of an inverse system $\{ \mathscr{F}_n \}_{n \in \mathbb{N}}$ such that
    \begin{enumerate}
        \item $\mathscr{F}_n$ is annihilated by $\mathscr{I}^{n+1}$, and
        \item The induced morphism $\mathscr{F}_{n+1} / \mathscr{I}^{n+1} \mathscr{F}_{n+1} \rightarrow \mathscr{F}_n$ is an isomorphism.
    \end{enumerate}
    With such description, one can show that the category of coherent $\mathscr{O}_X$-modules is Abelian and that exactness of sequences can be checked on an affine open cover (\cite[\href{https://stacks.math.columbia.edu/tag/087X}{Lemma 087X}]{stacks-project}).
\end{remark}

\begin{lemma}
\label{lem:local_exactness}
	Let $A$ be a Noetherian adic ring, let $M$ be an adic $A$-module and let $\mathfrak{m}$ be an open maximal ideal.
    Then $\hat{M}_{\mathfrak{m}}$, the $\mathfrak{m}$-adic completion of $M$, is isomorphic to $M \hat{\otimes}_A \hat{A}_{\mathfrak{m}}$, where $\hat{A}_{\mathfrak{m}}$ is the $\mathfrak{m}$-adic completion of $A$, and $M = 0$ if and only if $\hat{M}_{\mathfrak{m}} = 0$ for any open maximal ideal $\mathfrak{m} \subseteq A$.
    In particular, exactness of sequences of morphisms can be checked formal locally at open maximal ideals.
\end{lemma}

\begin{proof}
    The fact that $\hat{M}_{\mathfrak{m}} = M \hat{\otimes}_A \hat{A}_{\mathfrak{m}}$ follows from part (3) of \cite[\href{https://stacks.math.columbia.edu/tag/00MA}{Lemma 00MA}]{stacks-project}. \par

    Clearly, if $M = 0$, then $\hat{M}_{\mathfrak{m}} = 0$ for any open maximal ideal.
    Conversely, let $I$ be an ideal of definition of $A$, we will show that $M / I^{n+1} M = 0$ as a discrete $A / I^{n+1}$-module.
    Since $M$ is an adic $A$-module, this implies that $M = 0$.
    In order to show that $M / I^{n+1} M = 0$, it is enough to show that its $\mathfrak{m}$-adic completion is zero for all maximal ideals of $A / I^{n+1}$.
    Indeed, it certainly suffices to show that the localisation of $M / I^{n+1} M$ at every maximal ideal $\mathfrak{m}$ is zero and this can be checked on the $\mathfrak{m}$-adic completion of the local ring.
    Now consider the short exact sequence of $A$ modules
    \begin{align}
    \label{eq:short_exact_m_completion}
        0 \rightarrow I^{n+1} M \rightarrow M \rightarrow M / I^{n+1} M \rightarrow 0.
    \end{align}
    Taking the $\mathfrak{m}$-adic completion preserves exactness (part (2) of \cite[\href{https://stacks.math.columbia.edu/tag/00MB}{Lemma 00MB}]{stacks-project}), thus the $\mathfrak{m}$-adic completion of $M / I^{n+1} M$ is
    \begin{align}
    \label{eq:limit_quotient_commute}
        \frac{\hat{M}_{\mathfrak{m}}}{I^{n+1} \hat{M}_{\mathfrak{m}}}.
    \end{align}
    But maximal ideals of $A / I^{n+1}$ correspond bijectively to open maximal ideals of $A$ and (\ref{eq:limit_quotient_commute}) is $0$ by assumption, hence the claim. \par
    
    The last statement follows from the former as the category of adic $A$-modules is Abelian.
\end{proof}

\begin{lemma}
\label{lem:noetherian_closed_ideal}
	Let $X$ be a locally Noetherian formal scheme.
    Then every sheaf of ideals of $X$ is a coherent $\mathscr{O}_X$-module.
	In particular, there is a bijective correspondence between sheaves of ideals of $\mathscr{O}_X$ and closed formal subschemes of $X$.
\end{lemma}

\begin{proof}
	By \cite[Lemma \ref*{inf-lem:closed_immersion_local}]{bongiorno1}, it is sufficient to prove the correspondence locally.
    Hence let $X = \spf A$ and let $I \subseteq A$ be an ideal of definition.
    We only have to show that any ideal $K$ of $A$ is an adic $A$-module, the remaining statements will easily follow as in the case of schemes.
    Since $K$ is finitely generated, this amounts to showing that $K$ is a closed subspace of $A$, or equivalently that
	\begin{align}
	\label{eq:closed_condition}
		\bigcap_{n \in \mathbb{N}} I^{n+1} \left( A/K \right) = \emptyset.
	\end{align}
	By Krull's intersection theorem for Noetherian rings (part (2) of \cite[\href{https://stacks.math.columbia.edu/tag/00IQ}{Lemma 00IQ}]{stacks-project}), (\ref{eq:closed_condition}) holds if $I$ is contained in the Jacobson radical of $A$.
	Since $A$ is complete with respect to the $I$-adic topology, $I$ must be contained in every maximal ideal of $A$.
\end{proof}

\begin{example}
\label{ex:conormal_sheaf}
    Let $\iota : Z \subseteq X$ be a closed immersion of locally Noetherian formal schemes.
    The conormal sheaf $\mathscr{C}_{Z|X}$ is, by definition, the pull-back via $\iota$ of the ideal sheaf on $X$ defining $Z$.\
    By Lemma \ref{lem:noetherian_closed_ideal}, this is a coherent sheaf on $X$, hence $\mathscr{C}_{Z|X}$ is a coherent sheaf on $Z$.
    In particular, if $f : X \rightarrow Y$ is a morphism locally of formal finite presentation of locally Noetherian formal schemes, since $X \times_Y X$ is locally Noetherian, we see that $\Omega_{X/Y}^1$ is a coherent sheaf on $X$.
\end{example}

\subsection{Equivariant structures}
\label{subsec:equivariant_structures}

We follow \cite[\href{https://stacks.math.columbia.edu/tag/03LH}{Section 03LH}]{stacks-project} for basic notions of quasi-coherent sheaves on groupoids. \par

The following definition generalises the well-known notion of $G$-linearised vector bundles on $X$, for a fixed $G$-action on $X$, to the case of groupoids.

\begin{definition}
\label{def:linearisation}
    Let $\mathcal{X} = R \rightrightarrows X$ be an infinitesimal groupoid with both $X$ and $R$ locally Noetherian formal schemes.
    A \emph{coherent sheaf} on $\mathcal{X}$ is a coherent sheaf $\mathscr{F}$ on $X$ with an \emph{$R$-equivariant structure}, i.e. a choice of morphism
    \begin{align}
        \Phi : s^* \mathscr{F} \rightarrow t^* \mathscr{F}
    \end{align}
    of coherent sheaves on $R$ such that
    \begin{enumerate}
        \item $e^* \Phi : \mathscr{F} \rightarrow \mathscr{F}$ is the identity as quasi-coherent sheaves on $X$ (\emph{identity condition}), and
        \item the following diagram
        \begin{equation}
        \label{diag:cocycle}
            \begin{tikzcd}
                & \mathrm{pr}_2^* s^* \mathscr{F} \arrow[r, "\mathrm{pr}_2^* \Phi"] \arrow[dl, equals] & \mathrm{pr}_2^* t^* \mathscr{F} \arrow[dr, equals]& \\
                c^* s^* \mathscr{F} \arrow[dr, "c^* \Phi"'] & & & \mathrm{pr}_1^* s^* \mathscr{F} \arrow[dl, "\mathrm{pr}_1^* \Phi"] \\
                & c^* t^* \mathscr{F} \arrow[r, equals] & \mathrm{pr}_1^* t^* \mathscr{F} &
            \end{tikzcd}
        \end{equation}
        of coherent sheaves on $R \times_{(s,t)} R$ is commutative (\emph{cocycle condition}).
    \end{enumerate}
\end{definition}

If such an $R$-equivariant structure exists, $\Phi$ is necessarily an isomorphism (\cite[\href{https://stacks.math.columbia.edu/tag/077Q}{Lemma 077Q}]{stacks-project}).

\begin{example}
\label{ex:stack_differentials}
    Beside the trivial $R$-equivariant structure on $\mathscr{O}_X$, under smoothness hypotheses, there is another natural coherent sheaf associated to the groupoid with a natural $R$-equivariant structure.
    Suppose we are in Situation \hyperref[sit:global]{Global} and suppose further that $R \rightrightarrows X$ is an infinitesimal equivalence relation, i.e. the stabiliser group scheme is trivial over $X$.
    If we imagine the associated \emph{infinitesimal stack} $\mathcal{X} = [X / R]$ to be a genuine geometric space, we may wonder what its sheaf of differentials is, i.e. what is the conormal sheaf of the diagonal morphism $\mathcal{X} \rightarrow \mathcal{X} \times_S \mathcal{X}$.
    By working with a \emph{formally smooth presentation} $X \rightarrow \mathcal{X}$, we see that the following diagram
    \begin{equation}
    \label{diag:conormal_sheaf}
        \begin{tikzcd}[row sep=scriptsize, column sep=scriptsize]
            & R \arrow[dd, "\quad\quad s", sloped] \arrow[rr, "\chi"] \arrow[dl, "t", sloped] & & R \times_{(t,t)} R \times_{(s,t)} R \arrow[rr, "\mathrm{pr}_1 \times \mathrm{pr}_3"] \arrow[dd, "\quad\quad \Psi", sloped] \arrow[dl, "\mathrm{pr}_2", sloped] & & R \times_S R \arrow[dd, "s \times s", sloped]  \arrow[dl, "t \times t", sloped] \\
            X \arrow[rr, "\quad e", crossing over] & & R \arrow[rr, "\quad j", crossing over] & & X \times_S X & \\
            & X \arrow[rr, "e"]  & & R \arrow[rr, "j"] \arrow[dl] & & X \times_S X \arrow[dl] \\
            & & \mathcal{X} \arrow[rr, "\Delta_{\mathcal{X}/S}"] \arrow[from=uu, crossing over] & & \mathcal{X} \times_S \mathcal{X} \arrow[from=uu, crossing over] &
        \end{tikzcd}
    \end{equation}
    is commutative with $\chi(\gamma) = (\gamma, e(t(\gamma)), \gamma)$ and $\Psi(\alpha, \beta, \gamma) = \alpha^{-1} \circ \beta \circ \gamma$, and, furthermore, every face of the cube is Cartesian.
    Since $R$ is an infinitesimal equivalence relation, \cite[Lemma \ref*{frob-lem:monomorphism_closed}]{bongiorno2} implies that the induced morphism
    \begin{align}
    \label{eq:infinitesimal_diagonal}
        \hat{j} : R \hookrightarrow \left( X \times_S X \right)_{\hat{\Delta}},
    \end{align}
    where the codomain of $\hat{j}$ if the infinitesimal neighbourhood of the diagonal $\Delta : X \rightarrow X \times_S X$, is a closed immersion of formal schemes.
    Let $\mathscr{C}_{\hat{j}}$ be the conormal sheaf of $\hat{j}$ (Example \ref{ex:conormal_sheaf}), and define $\Omega_{\mathcal{X}/S}^1$ to be the coherent sheaf $e^*\mathscr{C}_{\hat{j}}$ on $X$ together with the $R$-equivariant structure hereafter described.
    By assumption, $s$ and $t$ are formally smooth, hence flat.
    Since conormal sheaves commute with flat base change (\cite[\href{https://stacks.math.columbia.edu/tag/0473}{Lemma 0473}]{stacks-project}) and infinitesimal neighbourhoods commute with base change (\cite[Lemma \ref*{inf-lem:infinitesimal_neighbourhood_fibre_products}]{bongiorno1}), we get the required natural isomorphism
    \begin{align}
        \label{eq:stack_differentials}
        s^* \Omega_{\mathcal{X}/S}^1 = s^* e^* \mathscr{C}_{\hat{j}} = \chi^* \mathscr{C}_{\hat{\mathrm{pr}}_1 \times \hat{\mathrm{pr}}_3} = t^* e^* \mathscr{C}_{\hat{j}} = t^* \Omega_{\mathcal{X}/S}^1.
    \end{align}
    An even larger diagram may be used to show that the cocycle condition is satisfied, however, we will refrain from elaborating further.
\end{example}

We move on to defining what the $R$-invariant subschemes of $X$ are.

\begin{definition}
\label{def:invariance}
    Let $\mathcal{X} = R \rightrightarrows X$ be an infinitesimal groupoid over a formal scheme $S$ and let $\iota : Z \rightarrow X$ be a formal subscheme of $X$.
    We say that $Z$ is \emph{$R$-invariant} if in the diagram
    \begin{equation}
    \label{diag:invariance}
        \begin{tikzcd}
            R|_Z \arrow[r, "\iota_t^{\prime}"] \arrow[d, "\iota_s^{\prime}"] & R \times_{(s, \iota)} Z \arrow[r, "s^{\prime}"] \arrow[d, "\iota_s"] & Z \arrow[d, "\iota"] \\
            Z \times_{(\iota, t)} R \arrow[r, "\iota_t"] \arrow[d, "t^{\prime}"] & R \arrow[r, "s"] \arrow[d, "t"'] & X \\
            Z \arrow[r, "\iota"] & X,
        \end{tikzcd}
    \end{equation}
    both (or equivalently either) $\iota_s^{\prime}$ and $\iota_t^{\prime}$ are isomorphisms.
\end{definition}

To check that the above definition makes sense, suppose that all of the above schemes are sets and that $R$ is an equivalence relation.
Then it is easy to see that $\iota_s^{\prime}$ is an isomorphism if and only if
\begin{align}
\label{eq:invariance_sets}
    x \sim_R y \text{ and } x \in Z \rightarrow y \in Z.
\end{align}
Similarly $\iota_t^{\prime}$ is an isomorphism if and only if (\ref{eq:invariance_sets}) holds after swapping $x$ and $y$ and clearly, using existence of $i$, if $\iota_s^{\prime}$ is an isomorphism, so is $\iota_t^{\prime}$.
Thus the notion of $R$-invariance is simply a scheme-theoretic upgrade.

\begin{remark}
\label{rem:invariance_properties}
    In general, restriction of a groupoid is not particularly well-behaved as the closed immersion $\iota : Z \rightarrow X$ is rarely formally smooth.
    Thus, if $R \rightrightarrows X$ is a formally smooth groupoid, the restriction $R|_Z \rightrightarrows Z$ need not be formally smooth.
    However, when $Z$ is $R$-invariant $R|_Z \rightrightarrows Z$ is formally smooth.
    More generally any property of a morphism stable by base change is preserved.
    This follows from observing Diagram (\ref{diag:invariance}) together with the fact that $\iota_s^{\prime}$ is an isomorphism.
\end{remark}

The next lemma is the first instance where we really do need the groupoid to be infinitesimal.

\begin{lemma}
    \label{lem:invariance_local}
    Let $R \rightrightarrows X$ be an infinitesimal groupoid over a formal scheme $S$ and let $\iota : Z \rightarrow X$ be a formal subscheme of $X$.
    Then the property of being $R$-invariant is affine-local on $X$.
\end{lemma}

\begin{proof}
    Certainly we may check whether $\iota_s$ is an isomorphism affine-locally on $R$.
    But since $e$ is a homeomorphism, $R$ is covered by open sets of the form $e(U)$ where $U$ is an affine open of $X$.
    Furthermore, since the image of an affine formal scheme under a thickening is affine, $e(U)$ is an affine formal subscheme of $R$.
\end{proof}

We can extend the notion of $R$-equivariance to morphisms $g : W \rightarrow X$ of formal schemes.

\begin{definition}
\label{def:equivariant_morphism}
    In Situation \hyperref[sit:global]{Global}, let $g : W \rightarrow X$ be a morphism locally of formal finite presentation of locally Noetherian formal schemes.
    An \emph{$R$-equivariant structure} on $g$, or on $W$, is a choice of morphism $\Phi$ preserving commutativity of the following diagram
    \begin{equation}
    \label{diag:equivariant_morphism}
        \begin{tikzcd}
            & R \times_{(s, g)} W \arrow[r, "s^{\prime}"] \arrow[d, "g_s"] \arrow[dl, "\Phi"'] & W \arrow[d, "g"] \\
            W \times_{(g, t)} R \arrow[r, "g_t"] \arrow[d, "t^{\prime}"] & R \arrow[r, "s"] \arrow[d, "t"'] & X \\
            W \arrow[r, "g"] & X,
        \end{tikzcd}
    \end{equation}
    which satisfies the corresponding identity and cocycle conditions of Definition \ref{def:linearisation}.
\end{definition}

\begin{remark}
\label{rem:equivariant_invariant}
    As in the case of coherent sheaves, $\Phi$ will necessarily be an isomorphism.
    Note that an $R$-equivariant structure on $W$ induces a commutative diagram
    \begin{equation}
    \label{diag:cartesian_equivariant}
        \begin{tikzcd}
            R_W^{\Phi} \arrow[r, shift left, "s_W"] \arrow[r, shift right, "t_W"'] \arrow[d, "g_{\Phi}"'] & W \arrow[d, "g"]\\
            R \arrow[r, shift left, "s"] \arrow[r, shift right, "t"'] & X,
        \end{tikzcd}
    \end{equation}
    where $W \times_{(g, t)} R =: R_W^{\Phi} := R \times_{(s, g)} W$, which is Cartesian with either pair of top or bottom arrows.
    We see that the groupoid $R_W^{\Phi} \rightrightarrows W$ inherits any property of the groupoid $R \rightrightarrows X$ which is stable under base change.
    Note how being $R$-invariant is a special case of being $R$-equivariant.
    Indeed, when $W = Z \subseteq X$ is a closed subscheme, there can be at most one equivariant structure on $Z$, or equivalently on $\mathscr{O}_Z$, or yet equivalently, by flatness of $s$ and $t$, on $\mathscr{I}_Z$.
    Such equivariant structure will necessarily satisfy the identity and cocycle conditions by uniqueness.
    In this case, we will simply say that $Z$, or $\mathscr{O}_Z$, or $\mathscr{I}_Z$, is $R$-invariant.
\end{remark}

The following are the main examples of equivariant structures we will need.
    
\begin{example}
\label{ex:push_presentation}
    Let $R \rightrightarrows X$ be an infinitesimal groupoid over a formal scheme $S$.
    Then the affine morphism $t : R \rightarrow X$ has a natural $R$-equivariant structure described by the following diagram
    \begin{equation}
    \label{diag:equivariant_push}
        \begin{tikzcd}
            & R \times_{(s, t)} R \arrow[r, "\mathrm{pr}_2"] \arrow[d, "\mathrm{pr}_1"] \arrow[dl, "c \times \mathrm{pr}_1"'] & R \arrow[d, "t"] \\
            R \times_{(t, t)} R \arrow[r, "\mathrm{pr}_2"] \arrow[d, "\mathrm{pr}_1"] & R \arrow[r, "s"] \arrow[d, "t"'] & X \\
            R \arrow[r, "t"] & X.
        \end{tikzcd}
    \end{equation}
    It is trivial to see that the diagram is commutative.
    In order to verify whether it satisfies the identity and cocycle conditions, we may pick an arbitrary $S$-scheme $T$ and reduce to verifying the case of sets.
    For instance, let us check that when base changing the commutative triangle in Diagram (\ref{diag:equivariant_push}) by $e : X \rightarrow R$, $\Phi$ restricts to $\mathds{1}_R$ (identity condition).
    For any $S$-scheme $T$, pick an element $(\alpha, \beta) \in \mathrm{Hom}_S(T, R \times_{(s,t)} R)$ and note that $\Phi(\alpha, \beta) = (\alpha \circ \beta, \alpha)$.
    After base changing, the morphism becomes $(e(t(\beta)), \beta) \rightarrow (\beta, e(t(\beta)))$, that is the identity morphism of $R$.
    Since we may test for all $S$-schemes, $\Phi$ indeed restricts to the identity.
    The cocycle condition may be checked similarly.
    We will refer to this $R$-equivariant structure on $t : R \rightarrow X$ as $R \times_{(s,t)} R \rightrightarrows R$ or simply as the natural $R$-equivariant structure. \par

    In fact, this equivariant structure can be thought as an equivariant structure on $t_* \mathscr{O}_R$, and, in stack-theoretic language, it may be described very concisely: if we let $\pi : X \rightarrow [X/R] =: \mathcal{X}$ denote the given formally smooth presentation, then the $R$-equivariant structure on $t_* \mathscr{O}_R$ is simply the quasi-coherent sheaf $\pi_* \mathscr{O}_X$ on $\mathcal{X}$.
\end{example}

\begin{example}
\label{ex:point_presentation}
    Now pick a closed point $x \in X$ and consider the formal subscheme $m : M := s^{-1}(x) \subseteq R$.
    We claim that the morphism $t : M \rightarrow X$ has a natural $R$-equivariant structure.
    Since $M$ is a formal subscheme of $R$ and $t : R \rightarrow X$ comes equipped with a natural $R$-equivariant structure, it suffices to check that $M$ is an $R \times_{(s,t)} R$-invariant subscheme of $R$.
    To this end, we pick an arbitrary $S$-scheme $T$ and we verify that $\Phi(T)$ restricts to a morphism of sets
    \begin{align}
    \label{eq:subset_invariant}
        \mathrm{Hom}_S \left( T, \left(R \times_{(s,t)} R \right) \times_{(\mathrm{pr}_2, m)} M \right) \rightarrow \mathrm{Hom}_S \left( T, M \times_{(m, \mathrm{pr}_1)} \left( R \times_{(t,t)} R \right) \right).
    \end{align}
    An element of the left-most set of (\ref{eq:subset_invariant}) can be described as a tuple $(\alpha, \beta)$ such that $s(\alpha) = t(\beta)$ and $s(\beta) = x$
    Under the morphism $\Phi(T)$, this is sent to $(\alpha \circ \beta, \alpha)$ and $s(\alpha \circ \beta) = s(\beta) = x$, hence $(\alpha \circ \beta, \alpha)$ is an element of the right-most set of (\ref{eq:subset_invariant}).
    It follows that $M$ is invariant, hence $t : M \rightarrow X$ has a natural $R$-equivariant structure. \par
    
    As in Example \ref{ex:push_presentation}, the natural $R$-equivariant structure on $M$ corresponds to a natural $R$-equivariant structure on $t_* \mathscr{O}_M$ and can be described as the quasi-coherent sheaf $\pi_* \mathscr{O}_x$ on $\mathcal{X}$, where $\pi : X \rightarrow [X/R] =: \mathcal{X}$ is the given formally smooth presentation.
\end{example}
\section{Functorial Constructions}
\label{sec:functorial_constructions}

In this section, we build several natural $R$-equivariant structures on $X$.
The general method is simple: suppose we have a functorial mapping within the category of formal schemes, for example normalisation or functorial resolution of singularities, which commutes with formally smooth base change, then there is a natural $R$-equivariant structure obtained by applying the functorial mapping to $R$.
Even if we are mainly interested in the case where $X$ is a usual scheme, we carry out all constructions in the category of formal schemes, as we need to be able to apply the mapping to the formal scheme $R$.

\subsection{Reduction}
\label{subsec:reduction}

Let $X$ be a locally Noetherian formal scheme.
We define the reduction of $X$ to be the closed formal subscheme $X_{\mathrm{red}} \subseteq X$ which is affine locally cut out by the nilradical ideal of $A$, where $X = \spf A$.
By Lemma \ref{lem:noetherian_closed_ideal}, $X_{\mathrm{red}}$ is a locally Noetherian formal scheme.
We remark that we are ignoring the topology of $X$ and that we are not quotienting by the topological nilradical.
It is not hard to see that the reduction satisfies the usual universal property and it is functorial.

\begin{lemma}
\label{lem:reduction_invariant}
    In Situation \hyperref[sit:global]{Global}, let $Z \subseteq X$ be an $R$-invariant closed formal subscheme of $X$.
    Then the reduction $Z_{\mathrm{red}} \subseteq X$ is also $R$-invariant.
\end{lemma}

\begin{proof}
    Up to restricting $R$ to the invariant subscheme $Z$ and renaming, it is sufficient to show that $X_{\mathrm{red}}$ is an $R$-invariant formal subscheme.
    Consider the following commutative diagram
    \begin{equation}
    \label{diag:reduction}
        \begin{tikzcd}
            R_{\mathrm{red}} \arrow[r] \arrow[d] & R \times_{(s, \iota)} X_{\mathrm{red}} \arrow[r] \arrow[d] & X_{\mathrm{red}} \arrow[d, "\iota"] \\
            X_{\mathrm{red}} \times_{(\iota, t)} R \arrow[r] \arrow[d] & R \arrow[r, "s"] \arrow[d, "t"'] & X \\
            X_{\mathrm{red}} \arrow[r, "\iota"] & X.
        \end{tikzcd}
    \end{equation}
    It suffices to show that
    \begin{align}
    \label{eq:reduction_isomorphism}
        R_{\mathrm{red}} \rightarrow R \times_{(s, \iota)} X_{\mathrm{red}}
    \end{align}
    is an isomorphism.
    In other words, we have to show that reduction commutes with formally smooth base change.
    We achieve this by reducing to the case of adic Noetherian local rings.
    Since the structure sheaves of both $R_{\mathrm{red}}$ and $R \times_{(s, \iota)} X_{\mathrm{red}}$ are coherent $\mathscr{O}_R$-modules, we may check that (\ref{eq:reduction_isomorphism}) is an isomorphism on an affine open cover of $R$.
    Thus we may assume that $R$, and hence $X$, are affine formal schemes.
    By Lemma \ref{lem:local_exactness}, we may choose a closed point $x \in R$, or equivalently a closed point $x \in X$, and show that (\ref{eq:reduction_isomorphism}) is an isomorphism after completing at the maximal ideals corresponding to $x$.
    This amounts to showing that reduction commutes with formally smooth base change for adic Noetherian local rings.
    This is achieved in Lemma \ref{lem:reduction_smooth_base_change}.
\end{proof}

\begin{lemma}
\label{lem:reduction_smooth_base_change}
    Let $\varphi : A \rightarrow B$ be a formally smooth local morphism of adic Noetherian local rings.
    Then if $A$ is reduced, so is $B$.
\end{lemma}

\begin{proof}
    We first observe that $\varphi$ is a regular morphism.
    This follows from the fact that $\varphi$ is a formally smooth morphism of adic Noetherian local rings together with the equivalences in \cite[\href{https://stacks.math.columbia.edu/tag/07PM}{Proposition 07PM}]{stacks-project}.
    Now, we apply Popescu's theorem (\cite[\href{https://stacks.math.columbia.edu/tag/07GB}{Section 07GB}]{stacks-project}) to express $\varphi$ as a filtered colimit of smooth (finitely presented) morphisms of rings $\varphi_{\lambda} : A \rightarrow B_{\lambda}$, for $\lambda$ in a filtered set $\Lambda$.
    By standard results for smooth morphisms, $B_{\lambda}$ is reduced (\cite[\href{https://stacks.math.columbia.edu/tag/033B}{Lemma 033B}]{stacks-project}) and it is easily verified that the colimit of reduced rings is itself reduced.
\end{proof}

\subsection{Normalisation}
\label{subsec:normalisation}

Let $X$ be a locally Noetherian formal scheme.
We define the normalisation of $X$ to be the morphism $\nu : X^{\nu} \rightarrow X$ affine locally given by the (discrete) normalisation $\nu^{\#} : A \rightarrow A^{\nu}$ (\cite[\href{https://stacks.math.columbia.edu/tag/035E}{Section 035E}]{stacks-project}), which we assume to be finite, where $X = \spf A$ and $A^{\nu}$ is endowed with the unique adic topology such that $\nu^{\#}$ is an adic morphism.
Under these assumptions, it is not hard to see that the normalisation satisfies the usual universal property and it is functorial.

\begin{lemma}
\label{lem:normalisation_equivariant}
    In Situation \hyperref[sit:global]{Global}, let $Z \rightarrow X$ be an $R$-invariant closed formal subscheme of $X$ whose normalisation $\nu : Z^{\nu} \rightarrow Z$ is finite.
    Then there exists a natural $R$-equivariant structure on $\nu : Z^{\nu} \rightarrow Z$.
\end{lemma}

\begin{proof}
    The proof is the same as in Lemma \ref{lem:reduction_invariant} after replacing $X_{\mathrm{red}}$ by $X^{\nu}$ and $\iota$ by $\nu$.
    We may still use Lemma \ref{lem:local_exactness} as, by assumption, $\nu_* \mathscr{O}_{X^{\nu}}$ is a coherent $\mathscr{O}_X$-module and we are reduced to show Lemma \ref{lem:normalisation_smooth_base_change}, which, in particular, implies that the normalisation of $R$ is finite.
    The identity and cocycle conditions are satisfied since normalisation and formally smooth base change commute.
\end{proof}

\begin{lemma}
\label{lem:normalisation_smooth_base_change}
    Let $\varphi : A \rightarrow B$ be a formally smooth local morphism of adic Noetherian local rings.
    Then if $A$ is normal, so is $B$.
\end{lemma}

\begin{proof}
    The proof is the same as in Lemma \ref{lem:reduction_smooth_base_change} but now we use the fact that, for a smooth morphism of rings $\varphi_{\lambda} : A \rightarrow B_{\lambda}$ with $A$ normal, $B_{\lambda}$ is normal (\cite[\href{https://stacks.math.columbia.edu/tag/033C}{Lemma 033C}]{stacks-project}), and that the filtered colimit of normal rings is normal.
    To see the latter assertion, first note that $B_{\lambda}$ is a direct sum of normal domains (\cite[\href{https://stacks.math.columbia.edu/tag/035Q}{Lemma 035Q}]{stacks-project}), thus, by commutativity of colimits, we may assume without loss of generality that each $B_{\lambda}$ is an integral domain.
    In this case, the colimit $B$ is also an integral domain and the assertion is an easy verification.
\end{proof}

\subsection{Fitting Ideals}
\label{subsec:fitting_ideals}

We would like to prove that the support of a sheaf with an equivariant structure is invariant.
However, when working with formal schemes, the set-theoretic support does not carry enough information, as any example with an adic Noetherian local ring confirms.
The scheme-theoretic upgrade of support is kindly provided by the theory of Fitting ideals.
We refer the reader to \cite[\href{https://stacks.math.columbia.edu/tag/07Z6}{Section 07Z6}]{stacks-project} and \cite[\href{https://stacks.math.columbia.edu/tag/0C3C}{Section 0C3C}]{stacks-project} for the construction of Fitting ideals in the case of schemes.

\begin{construction}
\label{cons:fitting_ideals}
    Let $X$ be a locally Noetherian formal scheme and let $\mathscr{F}$ be a coherent sheaf on $X$.
    We temporarily forget the topology of $\mathscr{F}$ and define its \emph{Fitting ideals} as the usual increasing sequence
    \begin{align}
    \label{eq:fitting_ideals}
        0 = \mathrm{Fitt}_{-1}(\mathscr{F}) \subseteq \mathrm{Fitt}_{0}(\mathscr{F}) \subseteq \mathrm{Fitt}_{1}(\mathscr{F}) \subseteq \ldots \subseteq \mathscr{O}_X.
    \end{align}
    By Lemma \ref{lem:noetherian_closed_ideal}, these correspond to a sequence of closed formal \emph{Fitting subschemes}
    \begin{align}
    \label{eq:fitting_subschemes}
        X = Z_{-1}(\mathscr{F}) \hookleftarrow Z_{0}(\mathscr{F}) \hookleftarrow Z_{1}(\mathscr{F}) \hookleftarrow \ldots \hookleftarrow \emptyset.
    \end{align}
    The results in \cite[\href{https://stacks.math.columbia.edu/tag/0C3C}{Section 0C3C}]{stacks-project} readily generalise to our context.
    In particular, we use the fact that formation of Fitting ideals commutes with flat base change (\cite[\href{https://stacks.math.columbia.edu/tag/0C3D}{Lemma 0C3D}]{stacks-project}) and the characterisation of locally free sheaves of constant rank through Fitting ideals (\cite[\href{https://stacks.math.columbia.edu/tag/0C3G}{Lemma 0C3G}]{stacks-project}).
\end{construction}

\begin{lemma}
\label{lem:fitting_invariant}
    In Situation \hyperref[sit:global]{Global}, let $\mathscr{F}$ be a coherent sheaf on $X$ endowed with an $R$-equivariant structure.
    Then the Fitting ideals of $\mathscr{F}$ are $R$-invariant.
\end{lemma}

\begin{proof}
    Let $\Phi$ denote the isomorphism $s^* \mathscr{F} \rightarrow t^* \mathscr{F}$ given by the $R$-equivariant structure.
    By commutativity of Fitting ideals and flat base change (\cite[\href{https://stacks.math.columbia.edu/tag/0C3D}{Lemma 0C3D}]{stacks-project}),
    \begin{align}
    \label{eq:fitting_ideals_extension}
        \mathrm{Fitt}_{r}(s^* \mathscr{F}) = s^* \mathrm{Fitt}_{r}(\mathscr{F}),
    \end{align}
    and similarly for $\mathrm{Fitt}_{r}(t^* \mathscr{F})$.
    On the other hand, since $\Phi$ is an isomorphism of sheaves over $R$, it induces an isomorphism of Fitting ideals
    \begin{align}
    \label{eq:fitting_isomorphism}
        \Phi_{\mathrm{Fitt}} : \mathrm{Fitt}_{r}(s^* \mathscr{F}) \xrightarrow{\sim} \mathrm{Fitt}_{r}(t^* \mathscr{F}).
    \end{align}
    This gives the required isomorphism between $s^* \mathrm{Fitt}_{r}(\mathscr{F})$ and $t^* \mathrm{Fitt}_{r}(t^* \mathscr{F})$.
\end{proof}

Despite the fact that $\Omega_{X/S}^1$ does not admit a natural equivariant structure, we can still show that its Fitting ideals are invariant.

\begin{lemma}
\label{lem:fitting_differentials_invariant}
    In Situation \hyperref[sit:global]{Global}, the Fitting ideals of $\Omega_{X/S}^1$ are $R$-invariant.
\end{lemma}

\begin{proof}
    Recall that there is a sequence of morphism of locally Noetherian formal schemes
    \begin{align}
    \label{eq:sequence_fitting}
        R \xrightarrow{s} X \rightarrow S,
    \end{align}
    where $s$ is formally smooth and locally of formal finite presentation.
    By \cite[Proposition 4.9]{MR2313672}, there exists a locally split short exact sequence of $\mathscr{O}_R$-modules
    \begin{align}
    \label{eq:short_sequence_fitting}
        0 \rightarrow s^* \Omega_{X/S}^1 \rightarrow \Omega_{R/S}^1 \rightarrow \Omega_{R/X}^1 \rightarrow 0,
    \end{align}
    where $\Omega_{R/X}^1$ is a locally free $\mathscr{O}_R$-module of finite rank $r$, for some $r \in \mathbb{N}$.
    Technically, the rank of $\Omega_{R/X}^1$ need not be constant, however we may test invariance of subschemes locally (Lemma \ref{lem:invariance_local}), hence we may assume $X$ connected and rank constant.
    We now fix $k \in \mathbb{N}$ and compute
    \begin{align}
    \label{eq:fitting_sum}
        \mathrm{Fitt}_{k}(\Omega_{R/S}^1) &= \sum_{i+j=k} \mathrm{Fitt}_{i}(s^* \Omega_{X/S}^1) \cdot \mathrm{Fitt}_{j}(\Omega_{R/X}^1) \nonumber \\
        &= \sum_{i=0}^{k-r} \mathrm{Fitt}_{i}(s^* \Omega_{X/S}^1) \nonumber \\
        &= \mathrm{Fitt}_{k-r}(s^* \Omega_{X/S}^1) \nonumber \\
        &= s^* \mathrm{Fitt}_{k-r}(\Omega_{X/S}^1),
    \end{align}
    where we employed part (2) of \cite[\href{https://stacks.math.columbia.edu/tag/07ZA}{Lemma 07ZA}]{stacks-project} and the local splitting property of (\ref{eq:short_sequence_fitting}) for the first equality, local freeness of $\Omega_{R/X}^1$ for the second equality, the sequence of inclusions in (\ref{eq:fitting_ideals}) for the third equality, and flatness of $s$ and invariance by base change of Fitting ideals for the fourth equality.    
    Repeating the argument with $t : R \rightarrow X$ and re-indexing (\ref{eq:fitting_sum}) yields the desired equality
    \begin{align}
    \label{eq:fitting_sum_conclusion}
        s^* \mathrm{Fitt}_{i}(\Omega_{X/S}^1) = \mathrm{Fitt}_{i+r}(\Omega_{R/S}^1) = t^* \mathrm{Fitt}_{i}(\Omega_{X/S}^1)
    \end{align}
    proving invariance of $\mathrm{Fitt}_{i}(\Omega_{X/S}^1)$, for any $i \in \mathbb{N}$.
\end{proof}

\subsection{Flattening Stratifications}
\label{subsec:flattening_stratifications}

Flattening stratifications are typically defined for morphisms of finite presentations, which is far from our situation.
Indeed, the flat locus of a local morphism of adic Noetherian local rings may fail to be open.
Nonetheless, there exists a local theory of flattening stratifications yielding a single maximal closed subscheme over which a given morphism is flat.
We follow \cite[\href{https://stacks.math.columbia.edu/tag/05LV}{Section 05LV}]{stacks-project} and we slightly simplify definitions to suit our case of interest: both rings are adic Noetherian local rings.

\begin{lemma}
\label{lem:flat_locus}
    Let $\varphi : A \rightarrow B$ be a local morphism of adic Noetherian local rings and let $\mathfrak{m}$ denote the maximal ideal of $A$.
    There exists a unique ideal $K_{\varphi} \subseteq \mathfrak{m}$ of $A$ such that
    \begin{enumerate}
        \item $B / K_{\varphi}B$ is flat over $A / K_{\varphi}$, and
        \item if $J \subseteq \mathfrak{m}$ satisfies (1), then $K_{\varphi} \subseteq J$.
    \end{enumerate}
\end{lemma}

\begin{proof}
    This is a special case of \cite[\href{https://stacks.math.columbia.edu/tag/0526}{Lemma 0526}]{stacks-project}.
\end{proof}

\begin{lemma}
\label{lem:flat_locus_base_change}
    Let $\varphi : A \rightarrow B$ be a local morphism of formal finite presentation of adic Noetherian local rings and let $\psi : A \rightarrow A^{\prime}$ be a formally smooth local morphism of formal finite presentation of adic Noetherian local rings.
    Assume that the completed fibre product $B^{\prime} := B \hat{\otimes}_A A^{\prime}$ is a local ring and let $\varphi^{\prime} : A^{\prime} \rightarrow B^{\prime}$ denote the base change morphism.
    Then $K_{\varphi^{\prime}} = \psi^* K_{\varphi}$.
\end{lemma}

\begin{proof}
    Since all morphisms are of formal finite presentation, the assumptions imply that $A$, $A^{\prime}$, $B$ and $B^{\prime}$ are adic Noetherian local rings and that $\psi$ is faithfully flat. \par

    We first show that $K_{\varphi^{\prime}} \subseteq \psi^* K_{\varphi}$, i.e. we show that the morphism $A^{\prime} / \psi^* K_{\varphi} \rightarrow B^{\prime} / \psi^* K_{\varphi} B^{\prime}$ is flat.
    To this end, we note that, by definition, $A / K_{\varphi} \rightarrow B/ K_{\varphi} B$ is flat and that, by Lemma \ref{lem:noetherian_closed_ideal}, its uncompleted base change is given by
    \begin{align}
    \label{eq:uncompleted_flat_locus}
        A^{\prime} / \psi^* K_{\varphi} \rightarrow \check{B} / \psi^* K_{\varphi} \check{B},
    \end{align}
    where $\check{B} := B \otimes_A A^{\prime}$.
    Certainly (\ref{eq:uncompleted_flat_locus}) is flat.
    But now $\check{B}$ is a Noetherian local ring, thus its completion $\check{B} / \psi^* K_{\varphi} \check{B} \rightarrow B^{\prime} / \psi^* K_{\varphi} B^{\prime}$ is faithfully flat and we get the desired containment. \par

    Conversely, we show that $K_{\varphi^{\prime}} \supseteq \psi^* K_{\varphi}$.
    Let $J$ be an ideal of $A^{\prime}$ such that $A^{\prime} / J \rightarrow B^{\prime} / J B^{\prime}$ is flat.
    We have to show that $\psi^* K_{\varphi} \subseteq J$, or equivalently, since $\psi$ is faithfully flat, $\psi(K_{\varphi}) \subseteq J$.
    To this end, we apply \cite[\href{https://stacks.math.columbia.edu/tag/0527}{Lemma 0527}]{stacks-project} to the Cartesian diagram
    \begin{equation}
    \label{diag:flat_locus_base_change}
        \begin{tikzcd}
            \check{B} / J \check{B} & B \arrow[l] \\
            A^{\prime} / J \arrow[u] & A \arrow[l, "\psi_J"] \arrow[u, "\varphi"]. 
        \end{tikzcd}
    \end{equation}
    The lemma concludes that $\psi_J(K_{\varphi}) = 0 \subseteq A^{\prime} / J$ so long as $\check{B} / J \check{B}$ is flat over $A^{\prime} / J$.
    This follows since the completion $B^{\prime} / J B^{\prime}$ is flat over $A^{\prime} / J$ by assumption, and is faithfully flat over $\check{B} / J \check{B}$, as $\check{B}$ is a Noetherian local ring.
\end{proof}

With the previous construction at hand, we define the maximal subscheme over which a given morphism is flat.

\begin{definition}
\label{def:flat_locus}
    Suppose that $g : W \rightarrow X$ is a morphism of formal finite presentation of Noetherian formal schemes whose topological spaces consist of a single closed point.
    Then we define the \emph{flat locus} $F_g$ of $g$ to be the non-empty formal subscheme cut out by the ideal constructed in Lemma \ref{lem:flat_locus}.
\end{definition}

\begin{lemma}
\label{lem:flat_locus_invariant}
    In Situation \hyperref[sit:local]{Local}, suppose that $g : W \rightarrow X$ is a morphism of formal finite presentation of Noetherian formal schemes whose topological spaces consist of a single closed point.
    Suppose further that $g$ is endowed with an $R$-equivariant structure $R_W^{\Phi}$ as in Diagram (\ref{diag:cartesian_equivariant}).
    Then $F_g$ is a non-empty $R$-invariant subscheme of $X$.
\end{lemma}

\begin{proof}
    Since $W$ is a Noetherian formal scheme whose topological space consists of a single closed point, so is $R_W^{\Phi}$, being an infinitesimal groupoid over $W$.
    Now, we can apply Lemma \ref{lem:flat_locus_base_change} to conclude that $s^* \mathscr{O}_{F_g} = \mathscr{O}_{F_{g_{\Phi}}} = t^* \mathscr{O}_{F_g}$.
\end{proof}

\subsection{Blowing up}
\label{subsec:blowing_up}

Let $X$ be a locally Noetherian formal scheme and let $\mathscr{K}$ be a sheaf of ideals on $X$.
We construct the blowing up of $X$ along $\mathscr{K}$.
This is done under very general finiteness conditions in \cite[\S 8.2, Definition 3]{MR3309387} but only for open ideals.
We prove that, under Noetherian assumptions, we can reduce to this construction. \par

We first construct the formal blowing up in the affine case.

\begin{construction}
\label{cons:blow_up_affine}
    Let $X = \spf A$, where $A$ is an adic Noetherian ring with ideal of definition $I$, and let $K$ be an ideal of $A$.
    We have already seen that $K$ is an adic $A$-module (Lemma \ref{lem:noetherian_closed_ideal}).
    We consider the Rees algebra $R_A(K) := \oplus_{m \in \mathbb{N}} \, K^m$ associated to $K$ and we endow it with the completion of the colimit topology of the subspace topologies on $K^m \subseteq A$ for $m \in \mathbb{N}$.
    Then, the natural morphism $A \rightarrow R_A(K)$ is an adic morphism of finite presentation and $R_A(K)$ is an adic Noetherian ring (Lemma \ref{lem:rees_topology}).
    This immediately implies that the ideal $I R_A(K)$ is an ideal of definition of $R_A(K)$. \par

    Let $X_n$ be the closed subscheme of $X$ cut out by $I^{n+1}$.
    We now define the blowing up of $X$ along the ideal $K$ to be $\tilde{X} = \lim{n \in \mathbb{N}} \tilde{X}_n$, where
    \begin{align}
    \label{eq:proj_blow_up}
        \pi_n : \tilde{X}_n := \mathrm{Proj}_{X_n} \left( R_A(K) \otimes_{A} A / I^{n+1} \right) \rightarrow X_n.
    \end{align}
    Note that $\tilde{X}_n$ is a well-defined scheme over $X_n$.
    Indeed, by the previous paragraph, $\tilde{X}_n$ is the $\mathrm{Proj}$ of a \emph{discrete} algebra over $A/I^{n+1}$ finitely generated in degree $1$.
    Now, the same construction as in \cite[\S 8.2]{MR3309387} applies, i.e. $\tilde{X}_n = X_n \times_{X_{n+1}} \tilde{X}_{n+1}$ and each $\tilde{X}_n \hookrightarrow \tilde{X}$ is a closed immersion for all $n \in \mathbb{N}$.
    It follows that $\tilde{X}$ is a Noetherian formal scheme and the sheaf of ideals associated to the closed immersion $\tilde{X}_0 \hookrightarrow \tilde{X}$ is a global ideal of definition.
    The limit of the morphisms $\pi_n$ is denoted by $\pi : \tilde{X} \rightarrow X$ and is the \emph{formal blowing-up} of $X$ along $K$.
\end{construction}

The following shows that, in a formal blowing up $\pi$, the inverse image of a closed (usual) subscheme of $X$ is a closed (usual) subscheme of $\tilde{X}$.

\begin{lemma}
\label{lem:rees_topology}
    Let $A$ be an adic Noetherian ring and let $K$ be an ideal of $A$.
    Then, the morphism $A \rightarrow R_A(K)$ to the Rees algebra of $K$ with its natural topology is an adic morphism of finite presentation and $R_A(K)$ is an adic Noetherian ring.
\end{lemma}

\begin{proof}
    We start by observing that $K$ is a closed finitely generated ideal of $A$ endowed with the $I$-adic topology (Lemma \ref{lem:noetherian_closed_ideal}), thus $A \rightarrow R_A(K)$ is finitely presented and $R_A(K)$ is Noetherian.
    Since $R_A(K)$ is, by construction, complete, to conclude, we have to show that $I R_A(K)$ is an ideal of definition of $R_A(K)$.
    Certainly, $I^{n+1} R_A(K)$ is open in $R_A(K)$, as its intersection with the degree $m$ component is $I^{n+1} \cap K^m$.
    This is, by definition of subspace topology, open, for any $n$ and $m$ in $\mathbb{N}$.
    Now let $J \subseteq R_A(K)$ be an open ideal.
    By definition $J_m := J \cap K^m$ is open for every $m \in \mathbb{N}$.
    Pick $m = 0$, then there exists an $n \in \mathbb{N}$ such that $I^{n+1} \subseteq J_0$.
    But then $I^{n+1} R_A(K) \subseteq J_0 R_A(K) \subseteq J$.
\end{proof}

Next, we show that the formal blowing up of an affine formal scheme possesses the well-known universal property.

\begin{lemma}
\label{lem:blow_up_affine_universal}
    Let $X = \spf A$, where $A$ is an adic Noetherian ring and let $K$ be an ideal of $A$.
    Let $\pi : \tilde{X} \rightarrow X$ denote the formal blowing up of $X$ along $K$.
    Then $\pi$ is universal (final) amongst all adic morphisms of finite presentation such that the inverse image of $K$ in $\tilde{X}$ is an invertible sheaf.
\end{lemma}

\begin{proof}
    This is the content of \cite[\S 8.2, Proposition 9]{MR3309387} in the affine case.
\end{proof}

We now finish off the construction in the global case.

\begin{construction}
\label{cons:blow_up}
    Let $X$ be a locally Noetherian formal scheme and let $\mathscr{K}$ be a sheaf of ideals.
    We construct the \emph{formal blowing up} of $X$ along $\mathscr{K}$ affine locally using Construction \ref{cons:blow_up_affine} and then we glue the resulting formal schemes using existence and uniqueness of glueing morphisms supplied by the universal property of the formal blowing up of affine formal schemes (Lemma \ref{lem:blow_up_affine_universal}).
    The result is an adic projective morphism $\pi : \tilde{X} \rightarrow X$ of finite presentation of locally Noetherian formal schemes.
\end{construction}

We record the universal property of formal blowing up in the more general non-affine case, however we will not use it.

\begin{lemma}
\label{lem:blow_up_universal}
    Let $X$ be a locally Noetherian formal scheme and let $\mathscr{K}$ be a sheaf of ideals on $X$.
    Let $\pi : \tilde{X} \rightarrow X$ denote the formal blowing up of $X$ along $\mathscr{K}$.
    Then $\pi$ is universal (final) amongst all adic morphisms of finite presentation such that the inverse image of $\mathscr{K}$ in $\tilde{X}$ is an invertible sheaf.
\end{lemma}

\begin{proof}
    This is the again the content of \cite[\S 8.2, Proposition 9]{MR3309387}.
\end{proof}

\begin{lemma}
\label{lem:blow_up_smooth_base_change}
    Let $X$ be a locally Noetherian formal scheme and let $\mathscr{K}$ be a sheaf of ideals of $X$.
    Suppose that $g : W \rightarrow X$ is a formally smooth morphism of formal finite presentation of locally Noetherian formal schemes and let $\pi : \tilde{X} \rightarrow X$ and $\pi^{\prime} : \tilde{W} \rightarrow W$ denote the formal blowing ups of $X$ and $W$ along $\mathscr{K}$ and $f^{-1} \mathscr{K}$ respectively.
    Then $\tilde{W} = W \times_X \tilde{X}$.
\end{lemma}

\begin{proof}
    Since $f$ is flat, we may directly apply \cite[\S 8.2, Proposition 4]{MR3309387} to conclude.
\end{proof}

\begin{lemma}
\label{lem:blow_up_invariance}
    In Situation \hyperref[sit:global]{Global}, let $\mathscr{K}$ be an $R$-invariant sheaf of ideals on $X$.
    Let $\pi : \tilde{X} \rightarrow X$ be the formal blowing up of $X$ along $\mathscr{K}$, then there exists a natural $R$-equivariant structure on $\pi$
    \begin{equation}
    \label{diag:blow_up_invariance}
        \begin{tikzcd}
            \tilde{R} \arrow[r, shift left, "\tilde{s}"] \arrow[r, shift right, "\tilde{t}"'] \arrow[d, "\tilde{\pi}"'] & \tilde{X} \arrow[d, "\pi"] \\
            R \arrow[r, shift left, "s"] \arrow[r, shift right, "t"'] & X,
        \end{tikzcd}
    \end{equation}
    where $\tilde{\pi} : \tilde{R} \rightarrow R$ is the blowing up of $R$ along $s^* \mathscr{K} = t^* \mathscr{K}$.
\end{lemma}

\begin{proof}
    This is immediate from the definition of invariant ideal and Lemma \ref{lem:blow_up_smooth_base_change}.
    The cocycle condition is satisfied since blowing up and formally smooth base change commute.
\end{proof}

\subsection{Functorial Resolutions}
\label{subsec:functorial_resolutions}

Recall from \cite[Theorem 1.2.1]{MR2957701}, that, for every reduced Noetherian quasi-excellent scheme $X$ over a field of characteristic zero, there exists a resolution $\tilde{X}$, which is functorial with respect to regular morphisms.
In more details, if $g : W \rightarrow X$ is a regular morphism of quasi-excellent schemes over a field of characteristic zero, then $\tilde{W} = W \times_X \tilde{X}$, where both $\tilde{X}$ and $\tilde{W}$ are regular.

\begin{lemma}
\label{lem:functorial_resolution}
    In Situation \hyperref[sit:local]{Local}, suppose further that $S = \spec k$ is a field of characteristic zero and that $X$ is reduced.
    Let $\pi : \tilde{X} \rightarrow X$ be the functorial resolution of $X$, then there exists a natural $R$-equivariant structure on $\pi$, as in Diagram (\ref{diag:blow_up_invariance}), where $\tilde{\pi} : \tilde{R} \rightarrow R$ is the functorial resolution of $R$.
\end{lemma}

\begin{proof}
    We only have to check that the results of \cite[Theorem 1.2.1]{MR2957701} apply.
    In our situation, $X$ and $R$ are adic Noetherian local rings, hence they are quasi-excellent; and $s$ and $t$ are formally smooth morphisms of adic Noetherian local rings, hence they are regular (\cite[\href{https://stacks.math.columbia.edu/tag/07PM}{Proposition 07PM}]{stacks-project}).
    The result now follows from functoriality.
\end{proof}
\section{Formal Leaves}
\label{sec:formal_leaves}

In this section, we define the formal leaf of an infinitesimal groupoid as the smallest invariant formal subscheme through a given closed point.
We then give a construction of formal leaf, which mimics the construction of the orbit of a groupoid, and we show that the absence of invariant subschemes forces pleasant regularity properties upon the formal leaf.
We conclude the section by giving a proof of existence of residual gerbes for algebraic stacks, which does not hinge on generic flatness.

\begin{definition}
\label{def:leaf}
    In Situation \hyperref[sit:local]{Local}, a \emph{leaf} of $R$ is a formal subscheme $l : L \hookrightarrow X$, which is $R$-invariant and does not contain any proper $R$-invariant formal subschemes.
    In Situation \hyperref[sit:global]{Global}, let $x \in X$ be a closed point.
    Then a \emph{leaf} of $R$ through $x \in X$ is defined by considering the infinitesimal neighbourhood of $x \in X$ and reducing to Situation \hyperref[sit:local]{Local}.
\end{definition}

\begin{lemma}
\label{lem:leaf_unique}
    In Situation \hyperref[sit:local]{Local}, the leaf of $R$ is unique.
\end{lemma}

\begin{proof}
    Let $L$ and $L^{\prime}$ be two leaves of $R$ through $x$ and consider their scheme-theoretic intersection $L \times_{X} L^{\prime}$.
    An easy categorical argument shows that $L \times_{X} L^{\prime}$ is also an $R$-invariant subscheme of $X$.
    Since it must contain $x$, it is equal to $L = L^{\prime}$.
\end{proof}

\begin{construction}
\label{cons:leaf}
    In Situation \hyperref[sit:local]{Local}, we construct the \emph{formal leaf} of $R$.
    Define the closed formal subscheme $l : L \hookrightarrow X$ to be the formal scheme-theoretic image of the morphism
    \begin{align}
    \label{eq:formal_leaf}
        t : M := s^{-1}(x) \rightarrow X
    \end{align}
    from Example \ref{ex:point_presentation}.
    Concretely, since $x$ is a closed point, $X = \spf A$ and $R = \spf \Gamma$ are affine formal schemes.
    Let $\mathfrak{m}$ be the maximal ideal of $A$ corresponding to $x$.
    Then $L$ is the closed formal subscheme corresponding to the ideal
    \begin{align}
        \label{eq:formal_leaf_ideal}
        K := \mathrm{ker} \left( A \xrightarrow{t^{\#}} \Gamma \otimes_{A} A /\mathfrak{m} \right) \subseteq A,
    \end{align}
    where the tensor product is via $s^{\#} : A \rightarrow \Gamma$.
\end{construction}

The next two lemmas show that $L$ is the formal leaf of $R$.

\begin{lemma}
\label{lem:leaf_invariant}
    In Situation \hyperref[sit:local]{Local}, the formal scheme $L$ of Construction \ref{cons:leaf} is an $R$-invariant formal subscheme of $X$.
\end{lemma}

\begin{proof}
    We consider Diagram (\ref{diag:equivariant_push}) after restricting to $M$ and inserting $L$:
    \begin{equation}
    \label{diag:leaf_invariant}
        \begin{tikzcd}
            & & R \times_{(s, t)} M \arrow[r, "\mathrm{pr}_2"] \arrow[d] \arrow[ddll, "c \times \mathrm{pr}_1"'] & M \arrow[d, "t"] \\
            & & R \times_{X} L \arrow[r, "\mathrm{pr}_2"] \arrow[d, "\mathrm{pr}_1"] \arrow[dl, dashed] & L \arrow[d, "l"] \\
            M \times_{(t, t)} R \arrow[r] \arrow[d, "\mathrm{pr}_1"] & L \times_{X} R \arrow[r, "\mathrm{pr}_2"] \arrow[d, "\mathrm{pr}_1"] & R \arrow[r, "s"] \arrow[d, "t"'] & X \\
            M \arrow[r, "t"] & L \arrow[r, "l"] & X. &
        \end{tikzcd}
    \end{equation}
    We have to show that the dashed morphism exists.
    To this end, it suffices to show that the formal scheme-theoretic image commutes with formally smooth, hence flat, base change.
    This is straightforward.
    Indeed, if the formal scheme-theoretic image of $M$ is cut out by an ideal $K$ inside $X$, the formal scheme-theoretic image of $R \times_{(s, t)} M$ is cut out by $s^* K$ inside $R$.
    This precisely describes the subscheme $R \times_{X} L$.
\end{proof}

\begin{lemma}
\label{lem:leaf_minimally_invariant}
    In Situation \hyperref[sit:local]{Local}, the formal scheme $L$ of Construction \ref{cons:leaf} contains no proper $R$-invariant formal subschemes.
\end{lemma}

\begin{proof}
    Suppose there exists an $R$-invariant formal subscheme $Z$ of $L$.
    Then $Z$ must contain the unique closed point $x$ and base changing the inclusion of $x$ in $Z$ by $s : R \rightarrow X$ yields a morphism $M \rightarrow R \times_{X} Z$.
    By assumption, $R \times_{X} Z = Z \times_{X} R$ as subschemes of $R$, therefore there exists a factorisation
    \begin{align}
    \label{eq:impossible_factorisation}
        t : M \rightarrow R \times_{X} Z = Z \times_{X} R \rightarrow Z \rightarrow X.
    \end{align}
    By definition of formal scheme-theoretic image of $t$, $Z = L$.
\end{proof}

We first prove some statements in the global setting, solely relying on the absence of proper $R$-invariant subschemes.

\begin{lemma}
\label{lem:coherent_equivariant_free}
    In Situation \hyperref[sit:global]{Global}, let $\mathscr{F}$ be a coherent sheaf on $X$ endowed with an $R$-equivariant structure.
    Suppose that $X$ has no proper $R$-invariant formal subschemes, then $\mathscr{F}$ is a locally free sheaf.
\end{lemma}

\begin{proof}
    By Lemma \ref{lem:fitting_invariant}, $\mathrm{Fitt}_{r}(\mathscr{F})$ are $R$-invariant sheaves of ideals.
    By assumption, either $\mathrm{Fitt}_{r}(\mathscr{F}) =0$ or $\mathrm{Fitt}_{r}(\mathscr{F}) = \mathscr{O}_X$ for any $r \in \mathbb{N}$.
    By \cite[\href{https://stacks.math.columbia.edu/tag/0C3G}{Lemma 0C3G}]{stacks-project}, $\mathscr{F}$ is locally free.
\end{proof}

\begin{lemma}
\label{lem:base_field}
    In Situation \hyperref[sit:global]{Global}, suppose that $X$ has no proper $R$-invariant formal subschemes.
    Then the image of the structural morphism $X \rightarrow S$ is the spectrum of the residue field of a closed point of $S$.
\end{lemma}

\begin{proof}
    We first show there exists a closed point $s \in S$ in the image of $X \rightarrow S$.
    We can easily reduce to the affine case by covering $S$ and $X$ by affine open subsets.
    Now, after quotienting by appropriate ideals of definition, the structural morphism becomes of finite presentation, and we may apply Chevalley's theorem (\cite[\href{https://stacks.math.columbia.edu/tag/00FE}{Theorem 00FE}]{stacks-project}) and the fact that $S$ is Noetherian to deduce existence of a closed point $s \in S$ in the image of $X \rightarrow S$.
    For such closed point $s \in S$ with residue field $\kappa(s)$, consider the base change $X_{\kappa(s)} \hookrightarrow X$.
    Since $R$ is a groupoid on $X$ over $S$, $X_{\kappa(s)}$ is an $R$-invariant subscheme of $X$, thus, by assumption, it must be equal to $X$.
    This proves the lemma.
\end{proof}

We finally prove regularity properties in the local setting.

\begin{lemma}
\label{lem:presentation_flat}
    In Situation \hyperref[sit:local]{Local}, the morphism $t : M \rightarrow L$ of Construction \ref{cons:leaf} is a faithfully flat morphism of formal finite presentation.
\end{lemma}

\begin{proof}
    Since $t : R \rightarrow X$ is a morphism of formal finite presentation, so is $t : M \rightarrow L$.
    Using Lemma \ref{lem:leaf_invariant} and extracting part of Diagram (\ref{diag:leaf_invariant}) yields existence of a commutative diagram
    \begin{equation}
    \label{diag:point_presentation}
        \begin{tikzcd}
            R_M \arrow[r, shift left, "s_M"] \arrow[r, shift right, "t_M"'] \arrow[d, "t^{\prime}"'] & M \arrow[d, "t"] \\
            R_L \arrow[r, shift left, "s_L"] \arrow[r, shift right, "t_L"'] & L,
        \end{tikzcd}
    \end{equation}
    which is Cartesian with both $s_L$ and $t_L$, where we renamed the groupoids $L \times_{X} R =: R_L := R \times_{X} L$ and $M \times_{(t, t)} R =: R_M := R \times_{(s, t)} M$.
    Now we apply Lemma \ref{lem:flat_locus_invariant} to conclude that the flat locus of $t : M \rightarrow L$ is an $R$-invariant formal subscheme of $L$, thus, by Lemma \ref{lem:leaf_minimally_invariant}, it has to be equal to $L$.
    Hence $t : M \rightarrow L$ is flat.
\end{proof}

\begin{lemma}
\label{lem:group_stabiliser}
    In Situation \hyperref[sit:local]{Local} and with the notation of Construction \ref{cons:leaf}, there is a Cartesian diagram
    \begin{equation}
    \label{diag:group_stabiliser}
        \begin{tikzcd}
            G \arrow[r, hook] \arrow[d] & M \arrow[d, "t"] \\
            \spec \kappa(x) \arrow[r, hook] & L,
        \end{tikzcd}
    \end{equation}
    where $G$ is the group stabiliser of $x \in X$.
\end{lemma}

\begin{proof}
    This follows from definition of group stabiliser and the fact that $M = R \times_X \spec \kappa(x)$.
\end{proof}

\begin{lemma}
\label{lem:leaf_regular}
    In Situation \hyperref[sit:local]{Local} and with the notation of Construction \ref{cons:leaf}, the formal leaf $L$ of $R$ is a regular formal scheme of dimension
    \begin{align}
    \label{eq:dimension_leaf}
        \mathrm{rk} \, R - \mathrm{dim} \, G.
    \end{align}
    Furthermore, the closed immersion $G \subseteq M$ of Diagram (\ref{diag:group_stabiliser}) is a local complete intersection.
\end{lemma}

\begin{proof}
    Note that $M$ is formally smooth over $\spec \kappa(x)$, hence it is regular (\cite[\href{https://stacks.math.columbia.edu/tag/07EI}{Lemma 07EI}]{stacks-project}), and, by definition, its dimension is $\mathrm{rk} \, R$.
    Since $t : M \rightarrow L$ is faithfully flat (Lemma \ref{lem:presentation_flat}), we see that $L$ is regular (\cite[\href{https://stacks.math.columbia.edu/tag/00OF}{Lemma 00OF}]{stacks-project}) of dimension $\mathrm{rk} \, R - \mathrm{dim} \, G$ (\cite[\href{https://stacks.math.columbia.edu/tag/00ON}{Lemma 00ON}]{stacks-project}), where we used the fact that the group stabiliser $G$ is the central fibre of $t$ (Lemma \ref{lem:group_stabiliser}).
    In particular, the closed immersion $\spec \kappa(x) \hookrightarrow L$ is a local complete intersection, and it follows that so must be its flat base change $G \hookrightarrow M$.
\end{proof}

\begin{lemma}
\label{lem:presentation_smooth}
    In Situation \hyperref[sit:local]{Local} and with the notation of Construction \ref{cons:leaf}, the morphism $t : M \rightarrow L$ is a formally smooth morphism of formal finite presentation if and only if $G$ is formally smooth over the residue field $\kappa(x)$.
\end{lemma}

\begin{proof}
    As already observed in Lemma \ref{lem:leaf_regular}, the central fibre of $t : M \rightarrow L$ is $G$.
    The lemma then follows from the fact that $t$ is faithfully flat (Lemma \ref{lem:presentation_flat}), so that $M$ is formally smooth over $L$ if and only if so is $G$ over $\kappa(x)$.
\end{proof}

We finally provide two proofs of existence of residual gerbes for Noetherian algebraic stacks.
Both proofs implicitly use the flattening stratification for adic Noetherian local rings.

\begin{lemma}
\label{lem:existence_gerbe}
    Let $\mathcal{X}$ be a Noetherian algebraic stack.
    Then the residual gerbe of $\mathcal{X}$ at $x$ exists for every $x \in | \mathcal{X} |$.
\end{lemma}

In the standard proof, one considers the morphism of stacks $\spec \kappa(x) \rightarrow \mathcal{X}$.
Up to replacing $\mathcal{X}$ with the stack-theoretic image, we assume the image of $\spec \kappa(x) \rightarrow \mathcal{X}$ is dense.
This is essentially the content of Construction \ref{cons:leaf} after choosing a smooth presentation.
In order to prove the lemma, we have to show that the morphism is flat.
Typically, this is accomplished invoking generic flatness.
See \cite[\href{https://stacks.math.columbia.edu/tag/0H22}{Lemma 0H22}]{stacks-project} for details.

\begin{proof}
[Proof 1]
    Since $\mathcal{X}$ is Noetherian, flatness can be checked on the completion.
    Base changing the morphism $\spec \kappa(x) \rightarrow \mathcal{X}$ by a smooth presentation yields the morphism $t : M \rightarrow L$ from Construction \ref{cons:leaf}.
    Now, Lemma \ref{lem:presentation_flat} shows flatness.
\end{proof}

Initially, we were trying to prove Lemma \ref{lem:presentation_flat} by expressing a quasi-coherent sheaf on a formal stack as a filtered colimit of coherent sheaves, thus reducing to Lemma \ref{lem:coherent_equivariant_free}.
It does not seem possible to always give such description in our generality, however, in the case of algebraic stacks, it is.

\begin{proof}
[Proof 2]
    Consider $\mathscr{O}_{\spec \kappa(x)}$ as a quasi-coherent sheaf over $\mathcal{X}$.
    We have to show it is flat.
    Since $\mathcal{X}$ is Noetherian, \cite[Proposition 15.4]{MR1771927} implies that $\mathscr{O}_{\spec \kappa(x)}$ can be expressed as $\mathrm{colim}\, \mathscr{F}_i$, where each $\mathscr{F}_i$ is a coherent sheaf on $\mathcal{X}$.
    By Lemma \ref{lem:coherent_equivariant_free}, each $\mathscr{F}_i$ is locally free, hence flat.
    The colimit of flat quasi-coherent sheaves is also flat (part (2) of \cite[\href{https://stacks.math.columbia.edu/tag/03EU}{Lemma 03EU}]{stacks-project}).
\end{proof}

Possibly, a similar strategy shows existence of residual gerbes for quasi-compact, quasi-separated algebraic stacks, after using \cite[Main Theorem]{MR3493431}.
\section{Minimal Presentations}
\label{sec:minimal_presentations}

In this section, we simplify the geometry of infinitesimal groupoids locally around a closed point by cutting the leaf with a transversal.
The infinitesimal groupoid restricted to the transversal will be locally equivalent to the original and, at the same time, it will have an invariant closed point.
We first define and construct the transversal and we then prove the main theorem, followed by auxiliary lemmas.

\begin{definition}
\label{def:transversal} 
    In Situation \hyperref[sit:global]{Global}, a \emph{transversal} of $R$ through a closed point $x \in X$ is a locally closed formal subscheme $v : V \rightarrow X$ containing $x$ intersecting the leaf $L$ of $R$ transversally locally around $x$.
    In more details, suppose that $\hat{V}$, the infinitesimal neighbourhood of $x \in V$, and that $\hat{V}$ and $L$ are cut out by ideals $\hat{J}$ and $K$ of $\hat{A}$, then
    \begin{enumerate}
        \setcounter{enumi}{-1}
        \item $\mathrm{Tor}_0^{\hat{A}} \left(\hat{A} / \hat{J}, \hat{A} / K \right) = \hat{A} / \hat{\mathfrak{m}}$, and
        \item $\mathrm{Tor}_1^{\hat{A}} \left(\hat{A} / \hat{J}, \hat{A} / K \right) = 0$,
    \end{enumerate}
    where $\hat{A}$ is the adic Noetherian local ring with maximal ideal $\hat{\mathfrak{m}}$ corresponding to the infinitesimal neighbourhood $\hat{X}$ of $x \in X$.
\end{definition}

\begin{remark}
\label{rem:transversal_regular_sequence}
    Condition (1) of Definition \ref{def:transversal} has been stated in such manner for simplicity and peace of mind, however, in practice, it will be more convenient to replace it by
    \begin{enumerate}
        \item[(1 bis)] there exist generators $\{ f_1, \ldots, f_c \}$ of $J$ whose images in $\hat{A}/K$, denoted by $\{\bar{f}_1, \ldots, \bar{f}_c \}$, form a regular sequence.
    \end{enumerate}
    It can be shown that (1) and (1 bis) are equivalent.
    Since we shall not use it, we do not prove it.
\end{remark}

\begin{construction}
\label{cons:transversal}
    It is straightforward to construct a transversal of $R$ through $x \in X$.
    Since we have to construct $V$ as a formal locally closed subscheme of $X$, we may assume $X = \spf A$ is affine.
    Let $\hat{A}$ denote the completion of $A$ at the maximal ideal $\mathfrak{m}$ corresponding to $x$ and note that, since $A$ is Noetherian, $\hat{A}$ is an adic Noetherian local ring with maximal ideal $\mathfrak{m} \hat{A} =: \hat{\mathfrak{m}}$ (\cite[\href{https://stacks.math.columbia.edu/tag/00MC}{Lemma 00MC}]{stacks-project}).
    Suppose further that the leaf $L$ is given by the formal spectrum of $\hat{A} / K$ and let $\mathfrak{n} := \mathfrak{m} \left( \hat{A} / K \right)$ be the maximal ideal corresponding to $x \in L$.
    By Lemma \ref{lem:leaf_regular}, $\mathfrak{n}$ is generated by a regular sequence, say $\{\bar{f}_1,...,\bar{f}_c\}$, over $\hat{A}/K$.
    Now lift $\{\bar{f}_1,...,\bar{f}_c\}$ to elements of $A$ (i.e. not $\hat{A}$), say $\{f_1,...,f_c\}$, and define $J$ to be the ideal of $A$ generated by such elements.
    This is possible since $\mathfrak{n}$ is generated by the image of $\mathfrak{m}$ in $\hat{A}/K$.
    In particular, we can ensure that $V := \spf A/J $ is a closed formal subscheme of $X$, and we write $\hat{V} = \spf \hat{A} / \hat{J}$ for the induced closed formal subscheme of $\hat{X}$.
    It remains to verify that 
    
    \begin{enumerate}[itemsep=1em]
        \setcounter{enumi}{-1}
        
        \item $\mathrm{Tor}_0^{\hat{A}} \left(\hat{A} / \hat{J}, \hat{A} / K \right) = \hat{A} / \hat{\mathfrak{m}}$, or equivalently $\hat{J} + K = \hat{\mathfrak{m}}$, which is true by construction, and

        \item $\mathrm{Tor}_1^{\hat{A}} \left(\hat{A} / \hat{J}, \hat{A} / K \right) = 0$, or equivalently (1 bis), which is true by construction and Lemma \ref{lem:leaf_regular}.
    \end{enumerate}
\end{construction}

Unlike the familiar non-singular case, being a transversal is not an open condition.

\begin{example}
\label{ex:transversal_not_open}
    Let $X = \mathbb{A}_k^2$ and let $R$ be the infinitesimal groupoid induced by the $\mathbb{G}_m$ action with weights $(1, 1)$.
    Then a transversal of $R$ through the origin is $V = \mathbb{A}_k^2$, and in fact, it is the only one.
    This follows from the fact that the origin is already $\mathbb{G}_m$-invariant.
    However, $V$ is not a transversal for any other point of $\mathbb{A}_k^2$, as any transversal would have to be $1$-dimensional.
\end{example}

\begin{theorem}
\label{thm:minimal_presentation}
    In Situation \hyperref[sit:global]{Global}, let $\mathcal{X}$ denote the formal stack $[X/R]$, let $x \in X$ be a closed point with residue field $\kappa(x)$ and suppose that its formal group stabiliser $G$ is formally smooth over $\kappa(x)$.
    Let $v : V \hookrightarrow X$ be a transversal of $R$ through $x \in X$, then, up to replacing $V$ with an open subset containing $x$,
    \begin{enumerate}[label=(\Alph*)]
        \item the closed point $x \in V$ is $R_V$-invariant, where $R_V$ is the restriction groupoid $R|_V$,
        \item the induced morphism $V \rightarrow \mathcal{X}$ is representable and formally smooth of formal finite presentation, and
        \item the immersion $v : V \hookrightarrow X$ is a local complete intersection.
    \end{enumerate}
    Suppose further that we are in Situation \hyperref[sit:local]{Local}, then
    \begin{enumerate}[resume, label=(\Alph*)]
        \item there exists a (non-unique) formally smooth morphism of formal finite presentation $u : X \rightarrow V$ of formal schemes over $\mathcal{X}$, split by $v$, and,
        \item if $G$ is trivial, $V = \mathcal{X}$ and $u : X \rightarrow V$ is unique.
    \end{enumerate}
    Suppose further that $S = \spec k$ is the spectrum of a field and that $x$ is a $k$-point.
    Let $l : L \hookrightarrow X$ be the leaf of $R$ through $x$.
    Then, after fixing $u : X \rightarrow V$,
    \begin{enumerate}[resume, label=(\Alph*)]
        \item there exists a (non-unique) decomposition of $X$ into a product $X \cong L \times_k V$ compatible with the morphisms $l : L \hookrightarrow X$, $v : V \hookrightarrow X$ and $u : X \rightarrow V$.
    \end{enumerate}
\end{theorem}

\begin{proof}
    We will first show part (A), (B) and (C) in Situation \hyperref[sit:local]{Local}, and denote the statements by (A'), (B') and (C') respectively.
    We will then prove the full statements in Situation \hyperref[sit:global]{Global}.
    For simplicity of notation, in parts (A'), (B') and (C'), we replace $X$, $V$ and $R$ with the infinitesimal neighbourhoods of $X$, $V$ and $R$ around $x$, $x$ and $e(x)$ respectively.
    We then consider the following commutative diagrams
    
    \begin{subequations}
        \begin{minipage}{0.55\textwidth}
        \centering
            \begin{equation}
            \label{diag:minimal_presentation}
                \begin{tikzcd}
                    G \arrow[r, hook] \arrow[d, hook] & M \arrow[r] \arrow[d, hook] & \spec \kappa(x) \arrow[d, hook] \\
                    N \arrow[r, hook] \arrow[d] & R \arrow[r, "s"] \arrow[d, "t"] & X \arrow[d] \\
                    V \arrow[r, "v", hook] & X \arrow[r] & \mathcal{X},
                \end{tikzcd}
            \end{equation}
        \end{minipage}
        \hfill
        \begin{minipage}{0.40\textwidth}
        \centering
            \begin{equation}
            \label{diag:minimal_presentation_zoom}
                \begin{tikzcd}
                    G \arrow[r, hook] \arrow[d, hook] & M \arrow[d, "t"] \\
                    \spec \kappa(x) \arrow[r, hook] \arrow[d, hook] & L \arrow[d, "l"] \\
                    V \arrow[r, "v", hook] & X,
                \end{tikzcd}
            \end{equation}
        \end{minipage}
    \end{subequations}
    where every square is Cartesian.
    Indeed, this is true for the rightmost diagram by condition (0) of Definition \ref{def:transversal} and Lemma \ref{lem:group_stabiliser}; for the leftmost diagram, $R = X \times_{\mathcal{X}} X$ and $M = R \times_X \spec \kappa(x)$ by definition, $N := V \times_X R$, and finally, using the rightmost diagram, we see that $G = V \times_X M$, which allows us to conclude that $G = N \times_R M$ via the pasting law for pull-backs. \par
    Furthermore, all formal schemes are formal spectra of adic Noetherian local rings and all morphisms are local morphisms.
    
    \begin{enumerate}[itemsep=1em, label=(\Alph*)]
        \item[(A')] Since $L$ is $R$-invariant, its intersection with $V$ is $R_V$-invariant.
        By Diagram (\ref{diag:minimal_presentation_zoom}), $V \cap L = \spec \kappa (x)$, showing that $x \in V$ is $R_V$-invariant.

        \item We have to show that, if $x \in \hat{V}$ is $\hat{R}_V$-invariant, then $x \in V$ is $R_V$-invariant.
        This is true since, for infinitesimal groupoids, being invariant is a local property (Lemma \ref{lem:invariance_local}).
        
        \item[(B')] From Diagram (\ref{diag:minimal_presentation}), we see that the base change of $V \rightarrow \mathcal{X}$ by the formally smooth presentation $X \rightarrow \mathcal{X}$ is the morphism $s : N \rightarrow X$.
        This is certainly a morphism of formal finite presentation of affine formal schemes, hence it suffices to show it is formally smooth.
        Since the central fibre of $s : N \rightarrow X$ is formally smooth by assumption, it suffices to show it is flat.
        To this end, we apply the slicing criterion of flatness (Lemma \ref{lem:slicing_criterion}) to the morphisms $s : N \hookrightarrow R \rightarrow X$.
        It suffices to find a sequence $\{f_1, \ldots, f_c \}$ of sections of $\mathscr{O}_R$ cutting out $N$, whose restriction to $M$ is a regular sequence.
        Take this sequence to be the extension via $t^{\#} : \mathscr{O}_X \rightarrow \mathscr{O}_R$ of the sequence described in Construction \ref{cons:transversal} to define $V$.
        Since the bottom left square of Diagram (\ref{diag:minimal_presentation}) is Cartesian, we see that $N \hookrightarrow R$ is cut out by such a sequence.
        Now, condition (0) of Definition \ref{def:transversal} implies that $\spec \kappa(x)$ is cut out by the restriction $\{\bar{f}_1,...,\bar{f}_c\}$ of the sequence $\{f_1, \ldots, f_c \}$ to $L$; and condition (1) of Definition \ref{def:transversal} implies that $\{\bar{f}_1,...,\bar{f}_c\}$ is a regular sequence.
        Finally, Diagram (\ref{diag:minimal_presentation_zoom}) implies that $G$ is cut out by the extension of $\{\bar{f}_1,...,\bar{f}_c\}$ via $t^{\#} : \mathscr{O}_L \rightarrow \mathscr{O}_M$ inside $M$.
        Since $t : M \rightarrow L$ is faithfully flat (Lemma \ref{lem:presentation_flat}), this is a regular sequence in $\mathscr{O}_M$ (\cite[\href{https://stacks.math.columbia.edu/tag/00LM}{Lemma 00LM}]{stacks-project}).
        Thus $s : N \rightarrow X$ is formally smooth.
        
        \item We know that $s : N \rightarrow X$ is a morphism of formal finite presentation whose completion at $x \in N$, $\hat{s}_x : \hat{N}_x \rightarrow \hat{X}_x$, is formally smooth.
        By \cite[Proposition 5.4]{MR2497584} together with the definition of \emph{smooth morphism} in \cite[1.3]{MR2497584}, there exists an open subset $U$ of $N$, and hence an open subset of $X$, containing $x$ over which $s$ is formally smooth.

        \item[(C')] By using Diagram (\ref{diag:minimal_presentation}) and faithfully flat descent for local complete intersections (\cite[\href{https://stacks.math.columbia.edu/tag/00LM}{Lemma 00LM}]{stacks-project}), we see that it suffices to show that $N \hookrightarrow R$ is a local complete intersection.
        By part (B'), this is a morphism of formal schemes over $X$ via $s$, where both $s : N \rightarrow X$ and $s : R \rightarrow X$ are formally smooth.
        Now Lemma \ref{lem:transversal_lci} implies that $N \hookrightarrow R$ is a local complete intersection.

        \item We have to show that, if $v : V \hookrightarrow X$ is a morphism of formal schemes whose completion at $x \in V$, $\hat{v}_x : \hat{V}_x \rightarrow \hat{X}_x$ is a local complete intersection, there exists an open subset $U$ of $V$ containing $x$ over which $v$ is a local complete intersection.
        This is true since, for Noetherian rings, being a regular ideal is a local property (\cite[\href{https://stacks.math.columbia.edu/tag/061L}{Lemma 061L}]{stacks-project}).
        
        \item We construct $u$ by first finding a dashed lift in the commutative diagram
        \begin{equation}
        \label{diag:presentation_lift}
            \begin{tikzcd}
                V \arrow[r, hook, "v"] \arrow[d, "\mathds{1}_V \times v"'] & X \arrow[d, "\mathds{1}_X"]  \arrow[ld, dashed] \\
                N \arrow[r, "s"] & X,
            \end{tikzcd}
        \end{equation}
        where we defined the morphism $V \rightarrow N = V \times_X R$ on each factor.
        Note that $v$ is a thickening of formal schemes and, in part (B'), we established that $s : N \rightarrow X$ is formally smooth.
        Therefore, the dashed morphism exists by \cite[Lemma \ref*{inf-lem:formally_smooth_lift_affine}]{bongiorno1}.
        We now define $u$ to be the composition of this dashed morphism with the projection $N \rightarrow V$.
        By construction, $v$ is a section of $u$ and, since $N$ is of formal finite presentation over $V$ and $X$ is a closed formal subscheme of $N$, $u$ is of formal finite presentation.
        Thus, it suffices to show that $u$ is formally smooth.
        By part (C'), $v$ is a local complete intersection, and, by construction, $v$ is a section of $u$, thus $u$ is flat (Lemma \ref{lem:lci_section_flat}).
        It follows that, in order to show that $u$ is formally smooth, it suffices to show that its central fibre is formally smooth.
        We consider the following commutative diagram
        \begin{equation}
        \label{diag:central_fibre_q}
            \begin{tikzcd}
                M \arrow[r, "t"] \arrow[d] & L \arrow[r, hook, "l"] \arrow[d] & X \arrow[d, "u"] \\
                G \arrow[r] \arrow[d] & \spec \kappa(x) \arrow[r] \arrow[d] & V \arrow[d] \\
                \spec \kappa(x) \arrow[r] & \left[ \ast / G \right] \arrow[r] & \mathcal{X},
            \end{tikzcd}
        \end{equation}
        where every square is Cartesian.
        This can be verified by either blindly trusting the expected properties of the residual gerbe or, slightly more rigorously, by noting that the formal scheme-theoretic image of $G \rightarrow V$ commutes with flat base change by $u$.
        Since $G$ is formally smooth, $t : M \rightarrow L$ is formally smooth and, by formally smooth descent (Lemma \ref{lem:formally_smooth_descent}), it follows that $L$ is formally smooth over $\kappa(x)$.
        
        \item If $G$ is trivial, then $R \rightrightarrows X$ is an infinitesimal equivalence relation (\cite[Section \ref{frob-sec:infinitesimal_equivalence_relations}]{bongiorno2}).
        The results then follow from part (1) and (4) of \cite[Main Theorem]{bongiorno2}.
        For more details, see Proposition \ref*{frob-prop:q_invariance}, Lemma \ref*{frob-lem:q_categorical} and Lemma \ref*{frob-lem:q_representable} from \cite{bongiorno2}.

        \item This is proved exactly as in part (5) and (6) of \cite[Main Theorem]{bongiorno2}.
        For more details, see Construction \ref*{frob-cons:alpha} and Lemma \ref*{frob-lem:alpha_isomorphism} from \cite{bongiorno2}.
    \end{enumerate}
\end{proof}

\begin{remark}
\label{rem:residual_gerbes}
    The above theorem implies that, in Situation \hyperref[sit:global]{Global}, the \emph{residual gerbe} exists.
    Indeed, by part (A), the closed point $x \in V$ is $R_V$-invariant, thus $[ \ast / G ]$ is a closed substack of $[V/R_V]$.
    By part (B), $V \rightarrow \mathcal{X}$ is a formally smooth presentation, hence $[V/R_V] = \mathcal{X}$ and $[ \ast / G ] \subseteq \mathcal{X}$.
\end{remark}

The remaining part of this section is devoted to proving the lemmas used in the proof of the theorem.

\begin{lemma}
\label{lem:slicing_criterion}
    Let
    \begin{equation}
    \label{diag:slicing_criterion}
        \begin{tikzcd}
            X \arrow[rr, hook, "\iota"] \arrow[rd, "p"'] & & Y \arrow[ld, "q"] \\
            & T &
        \end{tikzcd}
    \end{equation}
    be a commutative diagram of local morphisms of Noetherian formal schemes whose topological spaces consist of a single closed point.
    Let $t \in T$ be the unique closed point and let $f_t : X_t \hookrightarrow Y_t$ denote the base change of $f$ over $\spec \kappa(t)$.
    Assume that $\iota$ is a closed immersion and $q$ is flat.
    Suppose there exists a sequence of sections $\{f_1, \ldots f_c \}$ of $\mathscr{O}_Y$ cutting out $X$, whose restriction to $Y_t$ is a regular sequence.
    Then $p$ is flat.
\end{lemma}

\begin{proof}
    We already observed that all formal schemes are formal spectra of adic Noetherian local rings, hence we may directly apply \cite[\href{https://stacks.math.columbia.edu/tag/00MG}{Lemma 00MG}]{stacks-project}.
\end{proof}

\begin{lemma}
\label{lem:transversal_lci}
    Consider again Diagram (\ref{diag:slicing_criterion}), a commutative diagram of local morphisms of Noetherian formal schemes whose topological spaces consist of a single closed point.
    Assume that $\iota$ is a closed immersion and that $p$ and $q$ are formally smooth.
    Then $\iota$ is a local complete intersection.
\end{lemma}

\begin{proof}
    All formal schemes are formal spectra of adic Noetherian local rings, hence we may write Diagram (\ref{diag:slicing_criterion}) as $A \rightarrow B \rightarrow B/K =: C$ for an appropriate ideal $K$ of $B$.
    We now consider the fundamental triangle of cotangent complexes (\cite[\href{https://stacks.math.columbia.edu/tag/08QX}{Proposition 08QX}]{stacks-project}) associated to these morphisms
    \begin{align}
    \label{eq:fundamental_triangle}
        L_{B/A} \otimes_B^{\mathbf{L}} C \rightarrow L_{C/A} \rightarrow L_{C/B} \rightarrow L_{B/A} \otimes_B^{\mathbf{L}} C[1].
    \end{align}
    Since $p$ and $q$ are formally smooth morphisms of formal finite presentation, $L_{B/A}$ and $L_{C/A}$ are perfect complexes concentrated in degree $0$ (Lemma \ref{lem:formally_smooth_cotangent_complex_i}).
    But then it follow that $L_{C/B}$ is a perfect complex (\cite[\href{https://stacks.math.columbia.edu/tag/066R}{Lemma 066R}]{stacks-project}), which must be concentrated in degree $0$ and $-1$.
    Since it is a closed immersion, it is in fact perfect and concentrated in degree $-1$.
    Now, by the standard characterisation of local complete intersection morphisms for Noetherian rings (\cite[Corollary 6.14]{MR257068}), it follows that $K$ is a regular ideal of $B$.
\end{proof}

\begin{lemma}
\label{lem:lci_section_flat}
    Let $f : X \rightarrow Y$ be a local morphism of formal finite presentation of Noetherian formal schemes whose topological spaces consist of a single closed point.
    Suppose there exists a section $g : Y \rightarrow X$, which is a local complete intersection, then $f$ is flat.
\end{lemma}

\begin{proof}
    All formal schemes are formal spectra of adic Noetherian local rings, hence we may write $f$ as $A \rightarrow B$ and $g$ as $B \rightarrow B/K = A$ for an appropriate regular ideal $K$ of $B$.
    We see immediately that $B = A \oplus K$ as an $A$-module, hence it suffices to show that $K$ is a flat $A$-module.
    Note that the inverse limit of $K/K^{n+1}$ is flat.
    Indeed, $K$ is its completion with the respect to the $\mathrm{m}$-adic topology induced by $A$, and $\mathrm{m}$-adic completion of Noetherian rings preserves flatness.
    Since $K$ is generated by a regular sequence, $K^{n} / K^{n+1}$ is a flat $A$-module for any $n \in \mathbb{N}$, thus induction and \cite[\href{https://stacks.math.columbia.edu/tag/00HM}{Lemma 00HM}]{stacks-project} imply that $K / K^{n+1}$ is flat $A$-module for any $n \in \mathbb{N}$.
    But over a Noetherian ring, the inverse limit of flat modules with surjective transition morphisms is flat (\cite[\href{https://stacks.math.columbia.edu/tag/0912}{Lemma 0912}]{stacks-project}).
\end{proof}

\begin{lemma}
\label{lem:formally_smooth_descent}
    Consider a Cartesian diagram of local morphisms of formal finite presentation of Noetherian formal schemes whose topological spaces consist of a single closed point
    \begin{equation}
    \label{diag:formally_smooth_descent}
        \begin{tikzcd}
            X^{\prime} \arrow[d, "f^{\prime}"'] \arrow[r, "g^{\prime}"] & X \arrow[d, "f"] \\
            Y^{\prime} \arrow[r, "g"] & Y.
        \end{tikzcd}
    \end{equation}
    Assume that $g$ is flat.
    Then, if $f^{\prime}$ is formally smooth, so is $f$.
\end{lemma}

\begin{proof}
    By standard fpqc descent, we know that $f$ is flat, hence it suffices to show that its central fibre is formally smooth.
    As a result, after base changing to the unique closed point of $Y$, we may preserve all the assumptions of the lemma and assume that $Y = \spec k$ is a field.
    In this situation, by Lemma \ref{lem:formally_smooth_cotangent_complex_ii}, it is sufficient to show that the cotangent complex of $f$ is perfect and concentrated in degree $0$.
    By Lemma \ref{lem:formally_smooth_cotangent_complex_i}, the cotangent complex of $f^{\prime}$ is perfect and concentrated in degree $0$.
    Thus, by standard fpqc descent of perfect complexes (\cite[\href{https://stacks.math.columbia.edu/tag/09UG}{Lemma 09UG}]{stacks-project}) concentrated in degree $0$ (\cite[\href{https://stacks.math.columbia.edu/tag/09UD}{Lemma 09UD}]{stacks-project}), the cotangent complex of $f$ is perfect and concentrated in degree $0$.
\end{proof}

\begin{lemma}
\label{lem:formally_smooth_cotangent_complex_i}
    Let $\varphi : A \rightarrow B$ be a formally smooth local morphism of formal finite presentation of adic Noetherian local rings.
    Then its cotangent complex $L_{\varphi}$ is perfect and concentrated in degree $0$.
\end{lemma}

\begin{proof}
    If $\varphi$ is a formally smooth morphism, it is a regular morphism (\cite[\href{https://stacks.math.columbia.edu/tag/07PM}{Proposition 07PM}]{stacks-project}) and we may apply Popescu's theorem (\cite[\href{https://stacks.math.columbia.edu/tag/07GB}{Section 07GB}]{stacks-project}) to express $\varphi$ as a filtered colimit of smooth (finitely presented) morphisms of rings $\varphi_{\lambda} : A \rightarrow B_{\lambda}$, for $\lambda$ in a filtered set $\Lambda$.
    Certainly, the cotangent complex $L_{\varphi_{\lambda}}$ is concentrated in degree $0$ (\cite[\href{https://stacks.math.columbia.edu/tag/08R5}{Lemma 08R5}]{stacks-project}) for each $\lambda \in \Lambda$.
    By \cite[\href{https://stacks.math.columbia.edu/tag/08S9}{Lemma 08S9}]{stacks-project},
    \begin{align}
    \label{eq:cotangent_complex_colimit}
        L_{\varphi} = \colim{\lambda \in \Lambda} L_{\varphi_{\lambda}},
    \end{align}
    hence $L_{\varphi}$ is concentrated in degree $0$.
    It follows that $L_{\varphi}$ is quasi-isomorphic to $\Omega_{\varphi}^1$ and the latter sheaf is locally free.
\end{proof}

\begin{lemma}
\label{lem:formally_smooth_cotangent_complex_ii}
    Let $\varphi : A \rightarrow B$ be a local morphism of formal finite presentation of adic Noetherian local rings.
    Assume that $A = k$ is a field.
    Then, if $L_{\varphi}$ is perfect and concentrated in degree $0$, $\varphi$ is formally smooth.
\end{lemma}

\begin{proof}
    It suffices to show that $B$ is geometrically regular over $k$ (\cite[\href{https://stacks.math.columbia.edu/tag/07EL}{Theorem 07EL}]{stacks-project}).
    Let $\mathfrak{m}$ denote the maximal ideal of $B$ and let $K$ denote its residue field.
    Applying \cite[Corollary 5.10]{MR2497584} implies that the dimension of $B$ is equal to the dimension of $\mathfrak{m} / \mathfrak{m}^2$, thus $B$ is regular.
    In order to show it is geometrically regular, it suffices to show that the natural morphism $H^1(L_{K/k}) \rightarrow H^1(L_{K/A})  = \mathfrak{m} / \mathfrak{m}^2$ is injective (\cite[\href{https://stacks.math.columbia.edu/tag/07E5}{Proposition 07E5}]{stacks-project}).
    This follows by considering the exact sequence of cohomology of cotangent complexes associated to $k \rightarrow B \rightarrow K$, together with the assumption that $H^1(L_{A/k} \otimes_A^{\mathbf{L}} K) = 0$.
\end{proof}
\section{Foliations}
\label{sec:foliations}

In this section, we prove results about foliations.
In \S \ref{subsec:spaces}, we discuss the correspondence between infinitesimal groupoids and foliations.
In \S \ref{subsec:sheaves}, we discuss the correspondence between sheaves over such objects.
In \S \ref{subsec:foliation_singularities}, we discuss terminal foliation singularities.
In \S \ref{subsec:local_structure_theorems}, we prove the local structure theorems.

\begin{definition}
\label{def:foliation}
    Let $X$ be a locally Noetherian scheme over a field $k$ of characteristic zero.
    A \emph{foliation} $\mathscr{F}$ on $X$ over $k$ is a coherent subsheaf $\mathscr{F} \subseteq \mathscr{T}_{X/k}$ of the tangent sheaf which is \emph{involutive}, i.e. closed under the Lie bracket.
    If $X$ is a normal domain, we say $\mathscr{F}$ is \emph{saturated} if its normal sheaf $\mathscr{N}_{\mathscr{F}} := \mathscr{T}_{X/k} / \mathscr{F}$ is torsion-free.
\end{definition}

A foliation is therefore a Lie algebroid on $X$ over $k$ whose anchor map is injective.
For the notion of Lie algebroid and morphisms of Lie algebroids, we refer the reader to \cite[\S 6]{MR2157566}.

\subsection{Correspondence of Spaces}
\label{subsec:spaces}

The equivalence of categories between infinitesimal groupoids and Lie algebroids in characteristic zero is well-known to the experts.
Here below, we limit ourselves to collecting relevant results and verifying that the associated groupoid is formally smooth.

\begin{proposition}
\label{prop:foliation_groupoid}
    Let $X$ be a locally Noetherian scheme over a field $k$ of characteristic zero.
    Then there is an equivalence of categories
    \begin{align}
    \label{eq:foliation_groupoid}
        \left\{
            \begin{array}{c}
                R \rightrightarrows X \suchthat \text{$R$ is a formally smooth} \\
                \text{groupoid on $X$ over $k$ locally} \\
                \text{of formal finite presentation.}
            \end{array}
        \right\}
            \begin{array}{c}
                \xrightarrow{\;\alpha\;} \\
                \xleftarrow[\;\beta\;]{}
            \end{array}
        \left\{
            \begin{array}{c}
                \mathscr{F} \rightarrow \mathscr{T}_{X/k} \suchthat \text{$\mathscr{F}$ is a Lie} \\
                \text{algebroid on $X$ over $k$ and} \\
                \text{$\mathscr{F}$ is a locally free sheaf.}
            \end{array}
        \right\}
    \end{align}
\end{proposition}

As a sanity check, note that, when $X$ is a field of characteristic zero, we get the well-known correspondence between Lie algebras and formal groups.

\begin{proof}
    We describe the map $\alpha$.
    For a formally smooth groupoid of formal finite presentation $R \rightrightarrows X$, we consider the morphisms
    \begin{align}
        \label{eq:groupoid_to_foliation}
        \Delta_{X/k} : X \xrightarrow{e} R \xrightarrow{j} X \times_k X.
    \end{align}
    Since $e$ and $\Delta_{X/k}$ are closed immersions, we see that there is a morphism of conormal sheaves
    \begin{align}
        \label{eq:morphisms_conormal}
        \Omega_{X/k}^1 = \mathscr{C}_{\Delta_{X/k}} \rightarrow \mathscr{C}_e,
    \end{align}
    which is surjective if $j$ is a closed immersion (see \cite[\S 1 \& \S 2]{bongiorno2} for more details).
    Note that $\mathscr{C}_e = e^* \Omega_{s}^1$ is locally free, since $s$ is formally smooth.
    Define $\alpha(R \rightrightarrows X) = \mathscr{F}$ to be the dual of $\mathscr{C}_e$.
    Note that taking the dual of (\ref{eq:morphisms_conormal}) gives the anchor morphism $\mathscr{F} \rightarrow \mathscr{T}_X$.
    The exterior derivative naturally defined on $\Lambda^{\bullet} \, \mathscr{C}_e$, which is compatible with the exterior derivative of $\Lambda^{\bullet} \, \Omega_{X}$, gives rise to the Lie bracket on $\mathscr{F}$, which is compatible with the Lie bracket of $\mathscr{T}_X$. \par

    In fact, $\alpha$ induces an equivalence of categories.
    We first reduce to the case where $X$ is affine.
    Since $\mathscr{F}^{\vee}$ is a finite module over $X$, the source morphism $s$ is of formal finite presentation, hence $R$ is Noetherian.
    It follows that the property of being formally smooth of formal finite presentation is a local property (\cite[Corollary 5.5]{MR2497584}), thus it suffices to show formal smoothness affine locally on $X$.
    Finally, note that the global equivalence can be proved by glueing the local equivalences. \par
    
    In order to deal with the affine case, we let $R = \beta(\mathscr{F})$ be the infinitesimal groupoid constructed by applying \cite[Théorème 1.3.6, Chapitre VIII]{MR491681} to the dual of $\mathscr{F}$.
    This is allowed, for we are working in characteristic zero and the dual of a Lie algebroid has a natural compatible de Rham complex induced by its Lie bracket.
    Furthermore, we know that the identity morphism $e : X \rightarrow R$ is a quasi-regular immersion (\cite[Proposition 1.4.1.3, Chapitre VIII]{MR491681}), hence regular.
    It remains to verify that $R \rightrightarrows X$ is formally smooth of formal finite presentation.

    Since $X$ and $R$ are locally Noetherian, we may prove formal smoothness formal locally around closed points $x \in X$ (\cite[Corollary 5.5]{MR2497584}).
    Up to replacing $X$ and $R$, we may assume they are both Noetherian formal schemes whose topological spaces consist of a single closed point.
    Because $e$ is a regular immersion, we apply Lemma \ref{lem:lci_section_flat} to conclude that $s : R \rightarrow X$ is flat and it remains to show that its central fibre is formally smooth.
    To this end, we apply Lemma \ref{lem:functorial_foliations} to the natural equivariant structure on $t : R \rightarrow X$ from Example \ref{ex:push_presentation}, and, since $s$ is flat, we get an isomorphism of cotangent complexes $L_{\Delta_t} = s^* L_e$, where $\Delta_t$ is the diagonal morphism of $t : R \rightarrow X$.
    Using again flatness of $t$, it follows that $L_t$ is isomorphic to $L_{\Delta_t}$ shifted by $1$.
    Since $e$ is a regular immersion, $L_e$ is a perfect complex concentrated in degree $-1$ (\cite[Corollary 6.14]{MR257068}), thus $L_t$ is a perfect complex concentrated in degree $0$.
    Finally, since $t$ is flat and $L_t$ is perfect, \cite[\href{https://stacks.math.columbia.edu/tag/08QQ}{Lemma 08QQ}]{stacks-project} implies that the central fibre of $t$ has perfect cotangent complex concentrated in degree $0$.
    By Lemma \ref{lem:formally_smooth_cotangent_complex_ii}, we conclude that the central fibre of $t$ is formally smooth.
\end{proof}

\begin{lemma}
\label{lem:functorial_foliations}
    Let $R \rightrightarrows X$ be an infinitesimal groupoid on a scheme $X$ over a scheme $S$ and suppose that $g : W \rightarrow X$ is a morphism of schemes over $S$ endowed with an $R$-equivariant structure $Q \rightrightarrows W$.
    Let $\mathscr{C}_{e_X}$ and $\mathscr{C}_{e_W}$ be the conormal bundles of the identity morphisms $e_X : X \rightarrow R$ and $e_W : W \rightarrow Q$ respectively.
    Then there exists a natural isomorphism of Lie algebroids $g^* \mathscr{C}_{e_X} = \mathscr{C}_{e_W}$.
    Furthermore, if $g$ is flat, there is an isomorphism of cotangent complexes $g^* L_{e_X} = L_{e_W}$.
\end{lemma}

\begin{proof}
    In an ideal world, the proof should be an easy consequence of the fact that the commutative diagram
    \begin{equation}
    \label{diag:ideal_cartesian_stacks}
        \begin{tikzcd}
            W \arrow[r, "g"] \arrow[d, "\pi_W"'] & X \arrow[d, "\pi_X"] \\
            {[W/Q]} \arrow[r] & {[X/R]}
        \end{tikzcd}
    \end{equation}
    is Cartesian, together with the fact that $\Omega_{\pi_X}^1 = \mathscr{C}_{e_X}$.
    However, we will try to be slightly more rigorous.
    By definition of equivariant structure, we get the following commutative diagram
    \begin{equation}
    \label{diag:conormal_sheaves}
        \begin{tikzcd}
            W \arrow[r, "e_W"] \arrow[d, "g"] & Q \arrow[r, shift left, "s_W"] \arrow[r, shift right, "t_W"'] \arrow[d, "g"] & W \arrow[d, "g"] \\
            X \arrow[r, "e_X"] & R \arrow[r, shift left, "s_X"] \arrow[r, shift right, "t_X"'] & X,
        \end{tikzcd}
    \end{equation}
    where every square, with either $s$ or $t$, is Cartesian.
    We thus see that
    \begin{align}
    \label{eq:conormal_sheaves}
        \mathscr{C}_{e_W} = e_W^* \Omega_{s_W}^1 = e_W^* g^* \Omega_{s_X}^1 = g^* e_X^* \Omega_{s_X}^1 = g^* \mathscr{C}_{e_X}.
    \end{align}
    The statement about cotangent complexes is proved in the same way using flatness of $g$ and \cite[\href{https://stacks.math.columbia.edu/tag/08QQ}{Lemma 08QQ}]{stacks-project}.
\end{proof}

Now that we established how to move from one category to the other, we start to explore how to translate between properties of the foliation and properties of the groupoid.
We start with regularity.

\begin{definition}
\label{def:foliation_regular}
    Let $X$ be a locally Noetherian scheme over a field $k$ of characteristic zero and let $\mathscr{F}$ be a locally free foliation on $X$ over $k$.
    The \emph{regular locus} of $\mathscr{F}$ is defined to be
    \begin{align}
        \label{eq:regular_locus}
        \mathrm{Reg}(\mathscr{F}) = \left\{ x \in X \suchthat \Omega_{X,x}^1 \rightarrow \mathscr{F}_x^{\vee} \text{ is surjective.} \right\}
    \end{align}
    This is an open subset of $X$.
\end{definition}

It is well-known that, if $X$ is normal and $\mathscr{F}$ is saturated, $\mathrm{Reg}(\mathscr{F})$ contains all codimension $1$ points of $X$.

\begin{lemma}
\label{lem:stabiliser_surjection}
    Let $X$ be a locally Noetherian scheme over a field $k$ of characteristic zero and let $\mathscr{F}$ be a locally free foliation on $X$ over $k$.
    There is an exact sequence
    \begin{align}
    \label{eq:stabiliser}
        \Omega_{X/k}^1 \rightarrow \mathscr{F}^{\vee} \rightarrow \mathscr{C}_{X|P} \rightarrow 0,
    \end{align}
    where $\mathscr{C}_{X|P}$ is the conormal sheaf of the identity of the formal stabiliser group $P \rightarrow X$.
    In particular, $x \in \mathrm{Reg}(\mathscr{F})$ if and only if the formal stabiliser group of the associated infinitesimal groupoid $R \rightrightarrows X$ is trivial at $x$.
\end{lemma}

This lemma shows that the regular locus is precisely the locus where the diagonal morphism of the associated stack is unramified.
If furthermore $X$ is a normal domain and $\mathscr{F}$ is saturated, we deduce that all codimension $1$ points of $X$ have trivial formal group stabiliser.

\begin{proof}
    By definition, there is a commutative diagram
    \begin{equation}
    \label{diag:stabiliser}
        \begin{tikzcd}
            X \arrow[d, "\mathds{1}_X"'] \arrow[r, hook] & P \arrow[d] \arrow[r] & X \arrow[d, "\Delta_X"] \\
            X \arrow[r, "e", hook] & R \arrow[r, "j"] & X \times_k X,
        \end{tikzcd}
    \end{equation}
    where every square is Cartesian.
    It is easy to see this this implies existence of (\ref{eq:stabiliser}), where $\mathscr{C}_{X|P}$ is the conormal sheaf of the immersion $X \hookrightarrow P$.
    Hence $x \in \mathrm{Reg}(\mathscr{F})$ if and only if $\mathscr{C}_{X|P, x} = 0$ if and only if the stabiliser group $G_x = x$.
\end{proof}

\begin{example}
\label{ex:foliations}
We give a few examples of locally free foliations on locally Noetherian schemes over a field of characteristic zero.

    \begin{enumerate}[itemsep=1em]
        \item Applying (\ref{eq:foliation_groupoid}) to free Lie algebroids of rank $1$ with a choice of generator yields a correspondence
        \begin{align}
        \label{eq:vector_field_gm}
            \left\{
                \begin{array}{c}
                    \sigma : \hat{\mathbb{G}}_m \times X \rightarrow X \suchthat \text{$\sigma$ is a} \\
                    \text{formal group action.}
                \end{array}
            \right\}
                \begin{array}{c}
                    \xrightarrow{\;\alpha\;} \\
                    \xleftarrow[\;\beta\;]{}
                \end{array}
            \left\{
                \begin{array}{c}
                    \partial \in \mathscr{T}_{X/k} \suchthat \text{$\partial$ is a} \\
                    \text{vector field on $X$.}
                \end{array}
            \right\}
        \end{align}
        Indeed, given a vector field $\partial$, the free Lie algebroid of rank $1$ is obtained by letting $\mathfrak{g}_m$ be the abelian Lie algebra of dimension $1$ mapping to $\partial \in \mathscr{T}_X$, and tensoring $\mathfrak{g}_m \otimes_k \mathscr{O}_X \rightarrow \mathscr{T}_X$.
        Conversely, the identity immersion $X \rightarrow \mathbb{G}_m \times X$ has a distinguished first order deformation corresponding to an element of its normal bundle which naturally maps to $\mathscr{T}_X$ via $\sigma$ (see \ref{eq:morphisms_conormal}).
        This distinguished element is mapped to $\partial$.
        For formal groups, $\hat{\mathbb{G}}_m = \hat{\mathbb{G}}_a$, however, we have chosen to write $\hat{\mathbb{G}}_m$ as we want to think of vector fields as $1$-parameter subgroups of the foliation.
        
        \item Slightly more generally, given a Lie algebra $\mathfrak{g}$ acting on $X$, i.e. a morphism of Lie algebras $\mathfrak{g} \rightarrow \mathscr{T}_X$, we may construct the associated free Lie algebroid $\mathfrak{g} \otimes_k \mathscr{O}_X \rightarrow \mathscr{T}_X$.
        Via (\ref{eq:foliation_groupoid}), we think of a Lie algebra action as a formal group action $G \times X \rightarrow X$, where $G$ is the formal group associated to $\mathfrak{g}$.
        
        \item If $(X, D)$ is a log-smooth pair, $\mathscr{T}_X(-\mathrm{log} D) \subseteq \mathscr{T}_X$ is a locally free Lie algebroid.
        The induced foliation is almost never saturated, however, we may still construct its associated infinitesimal groupoid.
        Since a log-smooth pair is locally toric, we can think of this groupoid as \'{e}tale locally described by a formal $\hat{\mathbb{G}}_m^d$-action, where $d$ is the dimension of $X$.

        \item If the formal group stabiliser at a closed point $x \in X$ is semi-simple, the foliation is formal locally induced by a formal group action (\cite[Theoreme 2.2]{MR552968}), hence its double dual is locally free.
        
        \item If $X$ satisfies the hypotheses of the Zariski--Lipman conjecture, $\mathscr{T}_X$ is locally free and we can construct the associated infinitesimal groupoid.
    \end{enumerate}
\end{example}

\begin{remark}
\label{rem:foliation_resolution}
    Although there are plenty of locally free foliations around, most foliations arising from the Minimal Model Program are merely reflexive.
    There are two plausible routes to incorporate such foliations in the theory of stacks:

    \begin{enumerate}
        \item we can choose to work with infinitesimal \emph{derived} stacks.
        This is the solution employed in \cite{vezzosi}, or
        \item we can try to \emph{resolve} the singularities of the morphism $X \rightarrow [X/\mathscr{F}]$ and get an associated locally free foliation $\tilde{\mathscr{F}}$ on some birational model $\tilde{X}$, which, under good circumstances, should preserve the geometric properties of $\mathscr{F}$.
        This approach has been already carried out in some special cases in \cite{posva}, and in a differential geometric context in \cite{mohsen}.
    \end{enumerate}
    We will pursue (2) in future work.
\end{remark}

\subsection{Correspondence of Sheaves}
\label{subsec:sheaves}

We can generalise Proposition \ref{prop:foliation_groupoid} to a correspondence between coherent sheaves on infinitesimal groupoids and coherent sheaves with a flat partial connection valued in $\mathscr{F}$.

\begin{definition}
\label{def:partial_connection}
    Let $X$ be a locally Noetherian scheme over a field $k$ of characteristic zero and let $\mathscr{F}$ be a locally free Lie algebroid on $X$.
    Let $\mathscr{E}$ be a coherent $\mathscr{O}_X$-module.
    A \emph{flat partial connection} $\nabla$ on $\mathscr{E}$ valued in $\mathscr{F}$ is a morphism of sheaves over $k$
    \begin{align}
        \nabla : \mathscr{E} \rightarrow \mathscr{F}^{\vee} \otimes_{\mathscr{O}_X} \mathscr{E},
    \end{align}
    which, over any affine open subset $\spec A \subseteq X$, where $\mathscr{E}$ is given by $E$,
    \begin{enumerate}
        \item is $k$-linear, i.e. $\nabla(s + t) = \nabla(t) + \nabla(t)$ for all $s, t \in E$,
        \item satisfies Leibniz's rule, i.e. $\nabla(a s) = a \otimes \nabla(s) + da \otimes s$ for all $a \in A$ and $s \in E$; and
        \item is flat, i.e. $\nabla^2 = 0$, where the composition is uniquely defined.
    \end{enumerate}
\end{definition}

\begin{proposition}
\label{prop:equivariance_equivalence}
    Let $X$ be a locally Noetherian scheme over a field $k$ of characteristic zero and let $\mathscr{F}$ be a locally free Lie algebroid on $X$.
    Then there is an equivalence of categories
    \begin{align}
    \label{eq:connection_equivariance}
        \left\{
            \begin{array}{c}
                (\mathscr{E}, \Phi) \suchthat \text{$\mathscr{E}$ is a coherent} \\
                \text{sheaf on $X$ and $\Phi$ is an} \\
                \text{$R$-equivariant structure.}
            \end{array}
        \right\}
            \begin{array}{c}
                \xrightarrow{\;\alpha\;} \\
                \xleftarrow[\;\beta\;]{}
            \end{array}
        \left\{
            \begin{array}{c}
                \nabla : \mathscr{E} \rightarrow \mathscr{F}^{\vee} \otimes_{\mathscr{O}_X} \mathscr{E} \suchthat \text{$\mathscr{E}$ is a} \\
                \text{coherent sheaf on $X$ and} \\
                \text{$\nabla$ is a flat partial connection.}
            \end{array}
        \right\}
    \end{align}
\end{proposition}

\begin{proof}
    We have already seen such statement for $\mathscr{E} = \mathscr{O}_X$ with the foliated exterior derivative in the aforementioned \cite[Théorème 1.3.6, Chapitre VIII]{MR491681}.
    The same proof holds after tensoring by a coherent sheaf and utilising the given flat partial connection.
\end{proof}

Just as in the case of groupoids, we may consider what it means for a subscheme $Z \subseteq X$ to be $\mathscr{F}$-invariant.

\begin{definition}
\label{def:foliation_invariance}
    Let $X$ be a locally Noetherian scheme over a field $k$ of characteristic zero and let $\mathscr{F}$ be a locally free foliation on $X$.
    A subscheme $Z \subseteq X$, with ideal sheaf $\mathscr{I}_Z$, is $\mathscr{F}$-invariant if for all local derivations $\partial \in \mathscr{F}$, $\partial(\mathscr{I}_Z) \subseteq \mathscr{I}_Z$.
\end{definition}

\begin{lemma}
\label{lem:invariance_equivalence}
    Let $X$ be a locally Noetherian scheme over a field $k$ of characteristic zero and let $\mathscr{F}$ be a locally free foliation on $X$.
    Let $R \rightrightarrows X$ denote its associated infinitesimal groupoid.
    Then a subscheme $\iota : Z \hookrightarrow X$ is $\mathscr{F}$-invariant if and only if it is $R$-invariant.
\end{lemma}

\begin{proof}
    We first observe that $Z$ is $R$-invariant if and only if the dashed morphism in the commutative diagram
    \begin{equation}
    \label{diag:foliation_invariance}
        \begin{tikzcd}
            \iota^* \Omega_{X/k}^1 \arrow[d] \arrow[r, twoheadrightarrow] & \Omega_{Z/k}^1 \arrow[dl, dashed] \\
            \iota^* \mathscr{F}^{\vee}
        \end{tikzcd}
    \end{equation}
    exists.
    Indeed, the sufficient direction follows easily from Lemma \ref{lem:functorial_foliations} and the necessary direction follows from applying Proposition \ref{prop:foliation_groupoid} to the induced Lie algebroid $\iota^* \mathscr{F}$ on $Z$. \par

    Note that, if the dashed morphism exists it is unique.
    Hence, we localise the problem to $X = \spec A$, $Z = \spec A/I$ and $\mathscr{F}$ induced by the $A$-module $F$.
    Now the dashed morphism exists if and only if the exterior derivative $d : A \rightarrow F^{\vee}$ satisfies $d(I) \subseteq I F^{\vee}$.
    Finally, using the formula $\partial(x) = \partial(dx)$, where the latter $\partial$ is regarded as an element of $F^{\vee\vee} = F$, we conclude that $d(I) \subseteq I F^{\vee}$ if and only if $\partial(I) \subseteq I$.
\end{proof}

We can start applying some results proved in the previous sections to the case of foliations.

\begin{lemma}
\label{lem:sheaf_connection_free}
    Let $X$ be a locally Noetherian scheme over a field $k$ of characteristic zero and let $\mathscr{F}$ be a locally free foliation on $X$.
    Let $\mathscr{E}$ be a coherent sheaf with a flat partial connection.
    Then the Fitting ideals of $\mathscr{E}$ are $\mathscr{F}$-invariant and the restriction of $\mathscr{E}$ to the formal leaves of $X$ is locally free.
    In particular, a sheaf with a connection on a smooth scheme is locally free.
\end{lemma}

\begin{proof}
    After applying the correspondence of Proposition \ref{prop:equivariance_equivalence} and Lemma \ref{lem:invariance_equivalence}, the results follow from Lemma \ref{lem:fitting_invariant} and Lemma \ref{lem:coherent_equivariant_free}.
\end{proof}

\begin{lemma}
\label{lem:fitting_ideals_connection}
    Let $f : X \rightarrow Y$ be a morphism of locally Noetherian normal schemes of finite type over a field $k$ of characteristic zero and let $\partial$ be a local derivation of the relative tangent sheaf $\mathscr{T}_{X/Y}$.
    Let $\mathrm{Fitt}_{r}(\Omega_{X/Y}^1)$ be the $r\textsuperscript{th}$ Fitting ideal of $\Omega_{X/Y}^1$, then $\partial(\mathrm{Fitt}_{r}(\Omega_{X/Y}^1)) \subseteq \mathrm{Fitt}_{r}(\Omega_{X/Y}^1)$.
\end{lemma}

This result was proved in \cite[Corollary 2]{MR349654} for $Y = \spec k$.

\begin{proof}
    We apply the correspondence of Proposition \ref{prop:foliation_groupoid}, as in part (1) of Example \ref{ex:foliations}, so to get a formally smooth groupoid associated to $\partial \in \mathscr{T}_X$.
    Now the result follows from the characterisation of invariance of Lemma \ref{lem:invariance_equivalence} and the corresponding result for the Fitting ideals of the sheaf of differentials in Lemma \ref{lem:fitting_differentials_invariant}.
\end{proof}

\begin{example}
\label{ex:bott_connection}
We give a few examples of naturally arising coherent sheaves with a flat partial connection when the foliation $\mathscr{F}$ is locally free and $X$ is locally Noetherian over a field of characteristic zero.

    \begin{enumerate}[itemsep=1em]
        \item Let $D \subseteq X$ be an $\mathscr{F}$-invariant divisor.
        Then its ideal sheaf $\mathscr{I}_D$ is invariant, hence its dual $\mathscr{O}_X(D)$ has a natural flat partial connection.
        
        \item The conormal sheaf $\mathscr{C}_{X|P}$ from Lemma \ref{lem:stabiliser_surjection} describing the formal group stabiliser has a natural flat partial connection.
        This is induced from the sheaf of differentials of the diagonal of the associated stack $\Omega_j^1$, which was described in Example \ref{ex:stack_differentials}.
        Since we are only looking at the sheaf of differentials, as opposed to the degree $-1$ part of the cotangent complex, we obtain an equivariant structure even when the diagonal is not a closed immersion.
        Note that the Fitting ideals of $\mathscr{C}_{X|P}$ give a scheme-theoretic description of the non-regular locus of the foliation.
        By Lemma \ref{lem:sheaf_connection_free}, these are invariant by the foliation.
        On the other hand, the Fitting ideals of $\Omega_{X/k}^1$ give a scheme-theoretic description of the non-regular locus of the ambient scheme.
        By Lemma \ref{lem:fitting_ideals_connection}, these are also invariant by the foliation.

        \item We now obtain the \emph{Bott connection} as the flat partial connection associated to the sheaf of differentials of the associated stack.
        In what follows, $\mathscr{F}$ need not be locally free.
        Consider the open subset $U$ over which $X$ is regular and $\mathscr{F}$ is regular.
        When $X$ is normal and $\mathscr{F}$ is saturated, the closed complement of $U$ has codimension greater than or equal to $2$.
        When working over $U$, we apply Example \ref{ex:stack_differentials} to see that $e^* \mathscr{C}_j$ fits in a short exact sequence
        \begin{align}
            0 \rightarrow \mathscr{F} \rightarrow \mathscr{T}_X \rightarrow e^* \mathscr{C}_j \rightarrow 0.
        \end{align}
        of locally free sheaves over $U$.
        It follows that $e^* \mathscr{C}_j = \mathscr{N}_{\mathscr{F}}$ over $U$.
        Since $e^* \mathscr{C}_j$ has an equivariant structure, we get a flat partial connection on $\mathscr{N}_{\mathscr{F}}$ over $U$.

        \item When $X$ is normal and $\mathscr{F}$ is locally free and saturated, we get a flat partial connection on $\mathscr{C}_{\mathscr{F}} := \mathscr{N}_{\mathscr{F}}^{\vee}$, the \emph{conormal sheaf} of $\mathscr{F}$.
        This follows from observing that, by the previous part, there is flat partial connection
        \begin{align}
        \label{eq:flat_connection_conormal}
            \nabla_U : \mathscr{C}_{\mathscr{F}}|_U \rightarrow \mathscr{F}^{\vee}|_U \otimes \mathscr{C}_{\mathscr{F}}|_U
        \end{align}
        defined over $U$, and, since both sheaves are reflexive (we are crucially using that $\mathscr{F}^{\vee}$ is locally free), all sections of $\mathscr{F}^{\vee} \otimes \mathscr{C}_{\mathscr{F}}$ defined over $U$ extend to sections over $X$.
        Hence we can extend (\ref{eq:flat_connection_conormal}) to a flat partial connection on $X$.
        If $\mathscr{F}$ is not locally free, we still get a bilinear map
        \begin{align}
        \label{eq:flat_connection_normal}
            \mathscr{F} \otimes \mathscr{N}_{\mathscr{F}}^{\vee\vee} \rightarrow \mathscr{N}_{\mathscr{F}}^{\vee\vee}
        \end{align}
        by dualising (\ref{eq:flat_connection_conormal}) and extending the resulting section using the fact that $\mathscr{N}_{\mathscr{F}}^{\vee\vee}$ is reflexive.
        Since the connection is naturally defined, (\ref{eq:flat_connection_normal}) must be Bott connection.
    \end{enumerate}
\end{example}

\begin{remark}
\label{rem:bott_index_theorem}
    We show how, up to developing an intersection theory for infinitesimal stacks, the above correspondences give a convincing argument for \emph{Bott's vanishing formula}: let $X$ be a smooth projective variety of dimension $d$ over $\mathbb{C}$ with a regular foliation $\mathscr{F}$ of rank $r$ and let $D$ be an effective divisor such that $\mathscr{O}_X(D)$ admits a flat partial connection, then the self-intersection $D^k$ vanishes for all $k > d - r$. \par

    We give a sketch of a proof.
    Since $X$ is smooth and $\mathscr{F}$ is regular, $\mathscr{F}$ is locally free and we may apply the correspondences to get a formally smooth groupoid $R \rightrightarrows X$ and an $R$-equivariant line bundle $\mathscr{O}_X(D)$.
    By part (1) of Example \ref{ex:bott_connection}, this implies that $D$ is invariant.
    Let $\mathcal{X}$ denote the stack $[X/R]$.
    Since $D$ is $R$-invariant, $D$ descends to an effective divisor $\mathcal{D}$ on the stack.
    Since intersection products commute with flat pull-back, we may compute the $\mathcal{D}^k$ on $\mathcal{X}$.
    However, the dimension of $\mathcal{X}$ is $d -r$, hence $\mathcal{D}^k = 0$ for all $k > d - r$.
\end{remark}

\subsection{Foliation Singularities}
\label{subsec:foliation_singularities}

We will briefly make use of reflexive sheaves and the group of generalised divisors, as developed in \cite{MR2346188}.
We will still let $\otimes$ denote their group operation.

In the proof of Lemma \ref{lem:invariance_equivalence}, we saw that, if $\mathscr{F}$ is a locally free foliation on $X$ and $Z$ is an $\mathscr{F}$-invariant subscheme, then $\mathscr{F}|_Z$ is a locally free foliation on $Z$.
The next lemma shows that the same holds when blowing up, or, more generally, whenever we are given a morphism with an equivariant structure.

\begin{lemma}
\label{lem:blow_up_foliation}
    Let $X$ be a locally Noetherian scheme over a field $k$ of characteristic zero and let $\mathscr{F}$ be a locally free foliation on $X$.
    Let $Z$ be an $\mathscr{F}$-invariant subscheme of $X$ and let $\pi : \tilde{X} \rightarrow X$ be the normalisation of the blowing up of $X$ along $Z$.
    Then $\pi^* \mathscr{F}$ is a locally free (perhaps unsaturated) foliation on $\tilde{X}$.
\end{lemma}

When the rank of $\mathscr{F}$ is $1$, this was proved in \cite[Lemma I.1.3]{MR3644242}.

\begin{proof}
    We construct the infinitesimal groupoid $R \rightrightarrows X$ associated to $\mathscr{F}$ using Proposition \ref{prop:foliation_groupoid}.
    We saw that blowing up and normalising are functorial constructions (Lemma \ref{lem:blow_up_invariance} and Lemma \ref{lem:normalisation_equivariant} respectively), hence, since $Z$ is $R$-invariant (Lemma \ref{lem:invariance_equivalence}), $\pi$ is endowed with a natural equivariant structure $\tilde{R} \rightrightarrows \tilde{X}$.
    Now applying Lemma \ref{lem:functorial_foliations} and Proposition \ref{prop:foliation_groupoid} yields that $\pi^* \mathscr{F}$ is the foliation associated to the groupoid $\tilde{R}$ on $\tilde{X}$.
\end{proof}

The singularities of the Minimal Model Program measure the extent to which the induced foliation is unsaturated.
A more geometric way to see this is that they measure how far the induced formal group action on the exceptional divisor is from being trivial.
We illustrate this principle in the following lemma concerning foliated terminal singularities.
We refer the reader to \cite[\S 2.4]{MR4285142} for definition and properties of foliation singularities.

\begin{proposition}
\label{prop:terminal_singularities}
    Let $X$ be a normal domain of finite type over a field $k$ of characteristic zero and let $\mathscr{F}$ be a locally free saturated foliation on $X$.
    Suppose that $\mathscr{F}$ has at worst terminal singularities, then $X$ is regular and $\mathscr{F}$ is regular.
\end{proposition}

\begin{proof}
    Suppose that the dimension of $X$ is $d$ and that the rank of $\mathscr{F}$ is $r$. \par
    
    Suppose there exists an $\mathscr{F}$-invariant subscheme $Z \subseteq X$.
    We show that $Z$ must be a Cartier divisor inside $X$.
    Let $\pi : \tilde{X} \rightarrow X$ be the normalisation of the blowing up of $X$ along $Z$ and let $\tilde{\mathscr{F}}$ be the pull-back foliation on $\tilde{X}$.
    Applying Lemma \ref{lem:blow_up_foliation} to the pull-back construction yields that $\tilde{\mathscr{F}}$ is the saturation of $\pi^* \mathscr{F} \rightarrow \mathscr{T}_{\tilde{X}}$.
    Since $\tilde{X}$ is normal, we can discuss discrepancies.
    Using the fact that $\mathscr{F}$ is locally free, so that the dual operation commutes with tensor product, we see that the morphism $\pi^* \mathscr{F} \rightarrow \tilde{\mathscr{F}}$ induces a morphism
    \begin{align}
    \label{eq:canonical_singularities}
        \omega_{\tilde{\mathscr{F}}} = \left( \Lambda^r \tilde{\mathscr{F}} \right)^{\vee} \rightarrow \left( \Lambda^r \pi^* \mathscr{F}\right)^{\vee} = \pi^* \left( \Lambda^r \mathscr{F} \right)^{\vee} = \pi^* \omega_{\mathscr{F}}.
    \end{align}
    Since $\tilde{X}$ is normal, $\omega_{\tilde{\mathscr{F}}}$ is locally free, hence (\ref{eq:canonical_singularities}) is injective.
    But then
    \begin{align}
    \label{eq:terminal_singularities}
        \pi^* \omega_{\mathscr{F}} = \omega_{\tilde{\mathscr{F}}} \otimes \mathscr{O}_{\tilde{X}}(D),
    \end{align}
    for some divisor $D$ supported on the exceptional divisor $E$ of $\pi$.
    Since $\mathscr{F}$ is terminal, $E$ must be a divisor on $X$, i.e. $\pi$ is an isomorphism. \par

    Now, by Lemma \ref{lem:fitting_ideals_connection}, the Fitting ideals of $\Omega_{X/k}^1$ are $\mathscr{F}$-invariant.
    Since $X$ is normal, none of them can be locally principal.
    By the previous paragraph, they must all be $\mathscr{O}_X$.
    This shows that $\Omega_{X/k}^1$ is locally free (\cite[\href{https://stacks.math.columbia.edu/tag/0C3G}{Lemma 0C3G}]{stacks-project}), thus proving regularity of $X$. \par

    Finally, we consider the sheaf $\mathscr{C}_{X|P}$ from part (2) of Example \ref{ex:bott_connection}, together with its natural flat partial connection.
    By Lemma \ref{lem:sheaf_connection_free}, its Fitting ideals are $\mathscr{F}$-invariant.
    Since $\mathscr{F}$ is saturated, none of them can be locally principal.
    Applying the same argument as before shows that $\mathscr{C}_{X|P} = 0$, thus proving regularity of $\mathscr{F}$.
\end{proof}

\begin{proposition}
\label{prop:zariski_lipman}
    Let $X$ be a normal domain of finite type over a field $k$ of characteristic zero and suppose that its tangent sheaf $\mathscr{T}_{X/k}$ is locally free.
    Then, if $X$ has at worst terminal singularities, $X$ is regular.
\end{proposition}

This was conjectured in \cite{MR186672} for a general variety and proved in \cite[Theorem 6.1]{MR2854859} for Kawamata log terminal singularities, and in \cite[Theorem 1.1]{MR3239620} and \cite[Corollary 1.3]{MR3247804} for log canonical singularities.
We believe this proof is still worth writing since it does not make use of resolution of singularities.

\begin{proof}
    Clearly $\mathscr{T}_{X/k}$ is saturated and, by assumption, it is a foliation with at worst terminal singularities.
    We conclude by Proposition \ref{prop:terminal_singularities}.
\end{proof}

\subsection{Local Structure Theorems}
\label{subsec:local_structure_theorems}
We finally use Theorem \ref{thm:minimal_presentation} to deduce results about locally free foliations.

\begin{corollary}
\label{cor:smoothness_leaf}
    Let $X$ be a locally Noetherian scheme over a field $k$ of characteristic zero and let $\mathscr{F}$ be a locally free foliation on $X$.
    Then, for any closed point $x \in X$, the formal leaf of $\mathscr{F}$ through $x$ is regular.
\end{corollary}

\begin{proof}
    Apply Proposition \ref{prop:foliation_groupoid} and Lemma \ref{lem:leaf_regular}.
\end{proof}

\begin{corollary}
\label{cor:frobenius_foliation}
    Let $X$ be a locally Noetherian scheme over a field $k$ of characteristic zero and let $\mathscr{F}$ be a locally free foliation of rank $r$ on $X$.
    Suppose that, locally around a closed rational point $x \in X$, the induced morphism $\Omega_{X/k}^1 \rightarrow \mathscr{F}^{\vee}$ is surjective.
    Then there exists a locally closed subscheme $V \subseteq X$ containing $x$ such that, formal locally around $x \in X$, there is a decomposition
    \begin{align}
    \label{eq:frobenius_projection}
        X \cong \mathbb{A}_k^r \times V,
    \end{align}
    and $\mathscr{F}$ is the foliation induced by the projection $v : X \rightarrow V$.
\end{corollary}

This is a singular version of Frobenius theorem.

\begin{proof}
    Apply Proposition \ref{prop:foliation_groupoid} and part (E) and (F) of Theorem \ref{thm:minimal_presentation} on noting that surjectivity of $\Omega_{X/k}^1 \rightarrow \mathscr{F}^{\vee}$ in a neighbourhood of $x$ is equivalent to triviality of the stabiliser group at $x$ (Lemma \ref{lem:stabiliser_surjection}).
    The isomorphism $L \cong \mathbb{A}_k^l$ follows from Corollary \ref{cor:smoothness_leaf}.
\end{proof}

\begin{corollary}
\label{cor:zariski_foliation}
    Let $A$ be an adic Noetherian local $k$-algebra with maximal ideal $\mathfrak{m}$ and residue field $k$ of characteristic zero, and suppose there exists a derivation $\partial$ and an element $f \in \mathfrak{m}$ such that $\partial(f) \notin \mathfrak{m}$.
    Then there exists a subring $B \subseteq A$ and an isomorphism $A \cong B \llbracket x \rrbracket$ such that $\partial(B) = 0$.
\end{corollary}

This result was proved in \cite[Lemma 4]{MR177985}.

\begin{proof}
    We let $F$ be the sub-module of $\mathscr{T}_{A/k}$ generated by $\partial$.
    We show that $F$ is involutive, free, and the induced morphism $\Omega_{A/k}^1 \rightarrow F^{\vee}$ is surjective.
    Then we can directly apply Corollary \ref{cor:frobenius_foliation} to get the result. \par
    
    Clearly $F$ is involutive.
    To show it is free, we have to show that $a \partial \neq 0$ for any $0 \neq a \in A$.
    Write $u := \partial(f)$ for a unit $u$ in $A$.
    Then $a \partial(f) = au$, which is zero if and only if so is $a$.
    Finally, note that $df(\partial) = \partial(f) = u$, thus $u^{-1} dx$ is mapped to $\partial$ via $\Omega_{A/k}^1 \rightarrow F^{\vee}$.
\end{proof}

\begin{corollary}
\label{cor:cerveau_foliation}
    Let $X$ be a locally Noetherian scheme over a field $k$ of characteristic zero and let $\mathscr{F}$ be a locally free foliation on $X$.
    Then there exists a locally closed subscheme $V \subseteq X$ containing $x$ such that, formal locally around $x \in X$, there is a decomposition
    \begin{align}
    \label{eq:cerveau_projection}
        X \cong \mathbb{A}_k^l \times V,
    \end{align}
    where $l$ is the dimension of the formal leaf, $x \in V$ is $\mathscr{F}|_V$-invariant, $\mathscr{F}|_V$ is locally free, and $\mathscr{F}$ is the pull-back of $\mathscr{F}|_V$ via $v : X \rightarrow V$.
    If furthermore $X$ is a normal domain and $\mathscr{F}$ is saturated, so are $V$ and $\mathscr{F}|_V$.
\end{corollary}

This result was proved in \cite[Theoreme 1.1]{MR552968} for $X = \mathbb{A}_k^d$ and for foliations which are not necessarily locally free.

\begin{proof}
    Apply Proposition \ref{prop:foliation_groupoid} and all parts of Theorem \ref{thm:minimal_presentation} using the fact that, in characteristic zero, any formal group is formally smooth.
    The isomorphism $L \cong \mathbb{A}_k^l$ follows from Corollary \ref{cor:smoothness_leaf}. \par

    To check normality of $W$, we may work formal locally around $x \in W$ (\cite[\href{https://stacks.math.columbia.edu/tag/0FIZ}{Lemma 0FIZ}]{stacks-project}), where the claim follows from normality of $X$ and faithfully flat descent of normality (\cite[\href{https://stacks.math.columbia.edu/tag/033G}{Lemma 033G}]{stacks-project}) via $v$.
    To verify $\mathscr{F}|_V$ is saturated, we again work formal locally around $x$.
    Observe the morphism of short exact sequences
    \begin{equation}
    \label{diag:normal_sheaf}
        \begin{tikzcd}
            0 \arrow[r] & \mathscr{F} \arrow[d] \arrow[r] & \mathscr{T}_X \arrow[d, "dv"] \arrow[r] & \mathscr{N}_{\mathscr{F}} \arrow[d] \arrow[r] & 0 \\
            0 \arrow[r] & v^*\mathscr{F}|_V \arrow[r] & v^* \mathscr{T}_V \arrow[r] & v^* \mathscr{N}_{\mathscr{F}|_V} \arrow[r] & 0,
        \end{tikzcd}
    \end{equation}
    where the fact that $v$ is flat implies both that the lower sequence is exact and that $dv$ exists.
    Since $\mathscr{F}$ is the pull-back of $\mathscr{F}|_V$ via $v$, we have that $\mathscr{F} = dv^{-1}(v^*\mathscr{F}|_V)$, thus the morphism $\mathscr{N}_{\mathscr{F}} \rightarrow v^* \mathscr{N}_{\mathscr{F}|_V}$ is an isomorphism.
    Now, since $v$ is faithful and $\mathscr{N}_{\mathscr{F}}$ is torsion-free, we see that $\mathscr{N}_{\mathscr{F}|_V}$ must be torsion-free.
\end{proof}

\bibliographystyle{alpha}
\bibliography{main}

\end{document}